\documentclass[a4paper,10pt]{amsart}

\usepackage{geometry}

\usepackage[latin1]{inputenc}  
\usepackage[english]{babel}
\usepackage[T1]{fontenc} 
                                          
\usepackage{ams math}
\usepackage{amssymb}
\usepackage{mathrsfs}
\usepackage{amsfonts}
\usepackage{amsthm}
\usepackage{graphicx}
\usepackage{color}
\usepackage{soul}

\newtheorem{theo}{Theorem}[section]
\newtheorem{prop}[theo]{Proposition}

\newtheorem{lemme}[theo]{Lemma}
\newtheorem{corol}[theo]{Corollary}

\theoremstyle{remark}
\newtheorem{rem}[theo]{Remark}

\newcommand{\ve}{\varepsilon}
\newcommand{\vp}{\varphi}
\newcommand{\Vh}{\Vert_{H^1_h}}

\newcommand{\Lg}{\mathcal{L}_g}

\newcommand{\pui}{\frac{n-2}{2}}
\newcommand{\puiq}{\frac{2}{q_k-2}}

\newcommand{\RR}{\mathbb{R}}

\numberwithin{equation}{section}

\title{Effective multiplicity for the Einstein-scalar field Lichnerowicz equation}
\author{Bruno Premoselli}
\address{Laboratoire de Math\'ematiques - AGM, Universit\'e de Cergy-Pontoise, $2$ Avenue Adolphe Chauvin $95302$ Cergy-Pontoise Cedex and}
\address{Laboratoire de Math\'ematiques - UMPA, Ecole Normale Sup\'erieure de Lyon (site sciences), $46$ all\'ee d'Italie, $69364$ Lyon Cedex $07$. }

\email{bruno.premoselli@u-cergy.fr}

\begin{document}

\begin{abstract}
We prove the stability of the Einstein-scalar field Lichnerowicz equation under subcritical perturbations of the critical nonlinearity 
in dimensions $n = 3, 4, 5$. As a consequence, we obtain the existence of a second solution to the equation in several cases. In  particular, in the positive case, 
including the CMC positive cosmological constant case, we show that each time a solution exists, the equation produces a second solution with the exception of one critical 
value for which the solution is unique.
\end{abstract}

\maketitle

\section{Statement of the results.}

Let $(M,g)$ be a smooth closed Riemannian manifold of dimension $n \ge 3$. We are interested in the Einstein-scalar field Lichnerowicz equation in $M$:
\begin{equation} \label{EL} \tag{$EL$}
\triangle_g u + h u = f u^{2^*-1} + \frac{a}{u^{2^*+1}}\hskip.1cm ,
\end{equation}
where $h,f$ and $a$ are smooth functions on $M$, $\triangle_g = -\hbox{div}_g\nabla$ is the Laplace-Beltrami operator, $2^\star = \frac{2n}{n-2}$ 
is the critical Sobolev exponent for the Sobolev space $H^1$, and we assume that $\triangle_g + h $ is coercive, $a \ge 0$, $a \not \equiv 0$, and $\max_M f >0$. By closed, following standard terminology, we mean compact without boundary.
\medskip 

Equation \eqref{EL} arises in the mathematical analysis of  general relativity when solving the Einstein equations in a scalar-field setting, when the gravity is coupled to a scalar-field $\psi$.
Special important cases include the massive Klein-Gordon setting or the case of a positive cosmological constant $\Lambda$. Given a closed manifold $(M,g)$ of dimension $n \ge 3$ endowed with two smooth functions $\pi$ and $\psi$ and a $(2,0)$-symmetric tensor field $K$, the Cauchy problem in general relativity consists in finding a Lorentzian manifold $(M \times \RR,\tilde{g})$ together with a smooth function $\tilde{\psi}$ on $M \times \RR$ such that $\tilde \psi _{| M} = \psi$ and $\partial_n \tilde \psi_{|M} = \pi$, where $\partial_n$ denotes the normalized time derivative, such that $K$ is the second fundamental form of the embedding $M \subset M \times \RR$, and such that $(M \times \RR,\tilde{g})$ satisfies the Einstein equations:
\[ \hbox{Ric}(\tilde g)_{ab} + \frac{1}{2} R(\tilde g) \tilde{g}_{ab} = T_{ab}\hskip.1cm , \]
where $\hbox{Ric}(\tilde g)$ is the Ricci tensor of $\tilde g$, $R(\tilde g)$ is its scalar curvature and $T$ is the stress-energy tensor-field. This tensor field depends on $\tilde g$, on $\tilde \psi$ and on some potential $V$, itself related to $\tilde \psi$ by some wave equation. As shown first by Choquet-Bruhat \cite{ChoBru} for the vacuum case, see also Choquet-Bruhat-Isenberg-Pollack  \cite{ChoIsPo2}, a necessary and sufficient condition for the existence of such a $\tilde g$ on $M \times \RR$ is that the following system of equations in $M$ is satisfied:
\begin{equation} \label{sys:cont}
\begin{cases}
R(g) + \textrm{tr}_g K ^2  - \left| \left| K \right| \right|_{g}^2 & = \pi^2 + |\nabla \psi|_{g}^2 + 2V(\psi)\\
\partial_i (\mathrm{tr}_{g}K) - K_{i,j}^j  & = \pi \partial_i \psi \hskip.1cm ,\\
\end{cases} 
\end{equation}
where $R(g)$ is the scalar curvature of $g$ and $\nabla$ refers to the Levi-Civita connection of $g$. By specifying some of the unknown initial data $(g,K,\psi,\pi)$ and solving the system for the remaining data, the conformal method initiated by Lichnerowicz \cite{Lich} allows to turn \eqref{sys:cont} into a system of elliptic partial differential equations of critical Sobolev growth, called the conformal constraint system of equations. For a survey reference on the constraint equations see Bartnik-Isenberg \cite{BarIse} and for further informations on the conformal method see Choquet-Bruhat, Isenberg and Pollack \cite{ChoIsPo}. Essentially, the 
set of free data consists of $(\psi,\tau,\pi,U)$, where $\psi,\tau,\pi$ are smooth functions in $M$ and $U$ is a smooth symmetric traceless and divergence-free $(2,0)$-tensor in $M$. 
Given $( \psi, \tau, \pi, U)$ an initial free data set, the conformal constraint system 
of equations, whose unknowns are a smooth positive function $\vp$ in $M$ and a smooth vector field $W$ in $M$,  is written as
\begin{equation}\label{ELFull}
\begin{cases}
   \triangle_g \varphi + \mathcal{R}_\psi  \varphi  &= \mathcal{B}_{\tau, \psi, V} \varphi^{2^*-1} + \frac{\mathcal{A}_{\pi, U}(W)}{ \varphi^{2^*+1}}~,   \\ 
 \triangle_{g, conf} W  &= \frac{n-1}{n}\varphi^{2^*} \nabla\tau - \pi\nabla \psi\hskip.1cm , 
\end{cases}
\end{equation}
where we have let:
\begin{equation} \label{expressions}
\begin{aligned}
& \mathcal{R}_{\psi}  = c_n \left( R(g) - |\nabla \psi|_g^2 \right), \\
& \mathcal{B}_{\tau,\psi,V}  = c_n \left( 2 V(\psi) - \frac{n-1}{n} \tau^2 \right), \\
& \mathcal{A}_{\pi,U}(W)  = c_n \left(  |U + \mathcal{L}_g W |_g^2 + \pi^2 \right) ,\\
\end{aligned} \end{equation}
and $c_n = \frac{n-2}{4(n-1)}$. In the above $\triangle_{g, conf}W = \textrm{div}_g (\mathcal{L}_g W)$ and $\mathcal{L}_g W $ is the symmetric trace-free part of $\nabla W$ given by 
$ \mathcal{L}_gW_{ij} = W_{i,j} + W_{j,i} - \frac{2}{n} \mathrm{div}_g W g_{ij}$. Vector fields satisfying $\mathcal{L}_g W = 0$ are called conformal Killing vector fields. 
In the CMC case, which by definition corresponds to $\nabla\tau \equiv 0$, the system \eqref{ELFull} is semi-decoupled. The second equation in \eqref{ELFull} has a unique solution $W$ when $g$ has no conformal Killing vector fields. The first equation in \eqref{ELFull} is then nothing but \eqref{EL} with
$$h = \mathcal{R}_{\psi}~,~~f = \mathcal{B}_{\tau,\psi,V}~,~\hbox{and}~a = \mathcal{A}_{\pi,U}(W)~.$$
Note that conformal Killing vector fields generically do not exist by Beig, Chru\'sciel and Schoen \cite{BeChSc}. 
Several existence results for \eqref{EL} have been obtained. When $f \le 0$ the equation is fully understood, see Isenberg \cite{Ise} or Choquet-Bruhat, Isenberg and Pollack \cite{ChoIsPo}. 
Partial existence results are known when $\max_M f > 0$, see Hebey, Pacard and Pollack \cite{HePaPo} and Ng\^o and Xu \cite{NgoXu}. The main result of the paper is as follows.

\begin{theo} \label{Th1}
Let $(M,g)$ be a $n$-dimensional closed Riemannian manifold with $3 \le n \le 5$ and $h, f$ and $a$ be smooth functions on $M$. Assume that $\triangle_g+h$ is coercive, $a \ge 0$, $a \not \equiv 0$  and $\max_M f >0$. Consider the following Einstein-Lichnerowicz equation in $M$, for $\theta >0$:
\begin{equation} \label{ELt}\tag{$EL_\theta$}
\triangle_g u + hu = f u^{2^*-1} + \frac{\theta a }{u^{2^*+1}}.
\end{equation}
Then there exist $0 < \theta_1 \le \theta_2 \le + \infty$ such that equation \eqref{eq:einlicht} has:
\begin{itemize}
\item at least two solutions if $\theta < \theta_1$,
\item no solutions for $\theta > \theta_2$,
\item at least one solution if $\theta_1 \le \theta < \theta_2$.
\end{itemize}
In case $f >0$, $\theta_2$ is finite, there holds $\theta_1 = \theta_2$, and \eqref{ELt} has one and only solution for $\theta = \theta_1 = \theta_2$. In particular, when $f > 0$ in $M$, there exists $\theta_\star \in (0,+\infty)$ such that 
\eqref{eq:einlicht} has at least two solutions if $\theta < \theta_\star$, exactly one solution if $\theta = \theta_\star$, and no solution if $\theta > \theta_\star$.
\end{theo}

A few multiplicity results concerning \eqref{EL} were known. Ma and Wei \cite{MaWei} showed the existence of two smooth positive solutions of \eqref{EL} in dimensions $3 \le n \le 5$ 
when $h,f$ and $a$ are positive assuming the strict stability of the minimal solution but without investigating its stability.
In another context,  Ng\^o-Xu \cite{NgoXu} showed the existence of a second solution of \eqref{EL} if the first eigenvalue of the operator $\triangle_g + h$ is negative. Assuming that 
$\triangle_g+h$ is coercive, no results were known under the sole assumption $\max_Mf > 0$ that we use in Theorem \ref{Th1}. 
Assuming that $f$ is positive, Theorem \ref{Th1} shows that each time a solution exists, there is at least another solution to the equation, except in the limit-case $\theta = \theta_\star$. In particular, the $f > 0$ case in Theorem \ref{Th1} is proved by showing that the minimal solution of \eqref{EL} is always strictly stable, except in the limit-case $\theta = \theta_\star$, where it always fails to be strictly stable. The uniqueness in the limit case $\theta = \theta_\star$ is then obtained as a consequence of its non-strict stability.
This of course highly complements 
the Ma and Wei \cite{MaWei} result and shows that the Ma-Wei \cite{MaWei} alternative only occurs in the limit case $\theta = \theta_\star$. 
Needless to say, Theorem \ref{Th1} can be applied to the conformal constraint system \eqref{ELFull} in the CMC case. 
As a direct consequence 
of our theorem, noting that $\mathcal{A}_{\theta\pi,\theta U}(W) = \theta^2\mathcal{A}_{\pi,U}(W)$ in the CMC-case, we get that 
the following corollary holds true. 
The corollary includes, as a special case, the CMC positive cosmological constant case. 

\begin{corol} \label{CorolELSyst}
Let $(M,g)$ be a closed Riemannian manifold of dimension $3 \le n \le 5$ of positive Yamabe type such that $g$ has no conformal Killing vector fields. 
Let $V$ be a smooth positive function on $\RR$, and let  $\psi$ be a smooth function in $M$ such that the operator $\triangle_g + \mathcal{R}_\psi$ is coercive.
Assume that $(\tau, \pi, U)$ satisfies that $\tau \equiv C^{st}$ (CMC case), $(\pi,U) \not \equiv (0,0)$, and
\begin{equation} \label{signeB}
\frac{n-1}{n} \tau^2 < 2  \min_{x\in M}V \left( \psi(x) \right)~.
\end{equation} 
Then there exists $\theta_\star \in (0,+\infty)$ such that the conformal constraint system 
\begin{equation}\label{ELFullModif}
\begin{cases}
   \triangle_g \varphi + \mathcal{R}_\psi  \varphi  &= \mathcal{B}_{\tau, \psi, V} \varphi^{2^*-1} + \frac{\mathcal{A}_{\theta\pi,\theta U}(W)}{ \varphi^{2^*+1}}~,   \\ 
 \triangle_{g, conf} W  &= \frac{n-1}{n}\varphi^{2^*} \nabla\tau - \pi\nabla \psi~, 
\end{cases}
\end{equation}
has at least two solutions if $\theta < \theta_\star$, exactly one solution if $\theta = \theta_\star$, and no solution if $\theta > \theta_\star$.
\end{corol}

\begin{proof} Since $\mathcal{A}_{\theta \pi, \theta U}(W) = \theta^2 \mathcal{A}_{ \pi, U}(W)$ is nonnegative we first need to show that $\mathcal{A}_{ \pi, U}(W)$ is non zero, where $W$ is the unique solution of the vector equation in \eqref{ELFullModif} when $\nabla \tau = 0$. Using equation \eqref{expressions}, this is automatically true if $\pi \not \equiv 0$. If $\pi \equiv 0$ since $g$ has no conformal Killing fields this implies $W = 0$ and hence $U + \Lg W$ is not identically zero as soon as $U$ is not everywhere zero. 
Corollary \ref{CorolELSyst} then easily follows from Theorem \ref{Th1} since  \eqref{signeB} implies that $\mathcal{B}_{\tau,\psi,V} > 0$. 
\end{proof}

\medskip

The proof of Theorem \ref{Th1} goes through the proof of an involved stability result for \eqref{EL} under subcritical, asymptotically critical perturbations of the nonlinear power, and in particular makes use of 
blow-up analysis in the Sobolev critical setting. We state our stability result, Theorem \ref{Thstabi}, in section \ref{StateStabRes}. In section \ref{stabilite} we perform an asymptotic analysis of blowing-up sequences of solutions of 
our asymptotically critical equations to obtain sharp asymptotic estimates. They are used in section \ref{preuvethstabi} to prevent the appearence of concentration points and to prove Theorem \ref{Thstabi}.  In section \ref{mountainpass} we show that, if $3 \le n \le 5$, each time the critical equation \eqref{EL} 
has a mountain-pass structure it admits at least two smooth positive solutions. This result is obtained through a variational analysis of the subcritical equations obtained from \eqref{EL} and 
Theorem \ref{Thstabi}. Section \ref{solumini} is devoted to the construction of a minimal solution of \eqref{EL} for the $L^\infty(M)$-norm. 
Theorem \ref{Th1} is proved in section \ref{deuxsol}, using the results of sections \ref{mountainpass} and \ref{solumini}.

\subsection*{Acknowledgements.} 
The author warmly thanks Olivier Druet and Emmanuel Hebey for constant support and valuable remarks during the elaboration of this paper.

\section{Stability of the equation \eqref{EL}} \label{StateStabRes}

We present here the stability result we have for  \eqref{EL} which we will use to prove Theorem \ref{Th1}. Our stability result establishes that, in low dimensions, 
equation \eqref{EL} is stable under sub-critical perturbations of its critical exponent and perturbations of the coefficient $a$.

\begin{theo} \label{Thstabi}
Let $(M,g)$ be a $n$-dimensional closed Riemannian manifold with $3 \le n \le 5$ and $h, f$ and $a$ be smooth functions on $M$. Assume that $\triangle_g+h$ is coercive, $a$ is non-negative, $a \not \equiv 0$ and $\max_M f >0$. Let $(a_k)_k$ be a sequence of non-negative functions  converging in $C^0(M)$ to $a$ as $k \to \infty$ and 
$(q_k)_k$, $2 \le q_k \le 2^*$, be a sequence of positive real numbers converging to $2^*$. Consider the following (sub)critical perturbations of \eqref{EL}:
\begin{equation} \label{ELq} \tag{$EL_k$}
\triangle_g u + h u  = f u^{q_k-1} + \frac{a_k}{u^{q_k+1}}.
\end{equation}
Let $(u_k)_k$ be a sequence of solutions of \eqref{ELq} and assume that either $f > 0$ or $(u_k)_k$ is uniformly bounded in $H^1$. Then there exists a smooth positive solution $u$ of \eqref{EL} such that, up to a subsequence, $u_k$ converges to $u$ in the $C^{1, \alpha}$ topology for any  $0 < \alpha <  1$.
\end{theo}

Two stability results related to \eqref{EL} exist in the literature. Druet-Hebey proved in \cite{DruHeb}, in dimension $3$ to $5$, the stability and bounded-stability of \eqref{EL}
when perturbing the coefficients $h$, $f$ and $a$ without changing the critical exponent. In \cite{DruHeb} are also found examples of instability for \eqref{EL} when $n \ge 6$. Also Hebey-Veronelli \cite{HeVer} proved the stability of \eqref{EL} in the Einstein-Maxwell theory when $n = 3$ and $f = C^{st}$. The key point in Theorem \ref{Thstabi} with respect to \cite{DruHeb} and \cite{HeVer} is that we allow fully subcritical perturbations of \eqref{EL} which makes the analysis more involved as $q_k \to 2^*$ with $q_k < 2^*$ for all $k$. This is crucial in order to obtain the existence of two solutions in Theorem \ref{Th1}, which are constructed as limits of solutions of the subcritical equations associated to \eqref{EL}, see section \ref{mountainpass}.

\medskip
As an interesting remark it can be noted that the generic existence of two solutions in Theorem \ref{Th1} is implicitly contained in Theorem \ref{Thstabi} through a degree-theory argument. 
Consider for instance the case $f > 0$. We compute the degree of \eqref{ELt} mimicking a computation carried out by Schoen \cite{SchoenNumberMetrics} 
for the Yamabe equation. Since $f >0$ both the existence of a 
solution of \eqref{ELt} when $\theta \ll 1$ and the non-existence when $\theta \gg 1$ were proven in Hebey-Pacard-Pollack \cite{HePaPo}. Let for any positive $M$ 
\begin{equation*}
\Omega_M = \left \{ u \in C^2(M) \hbox{ such that } \| u\|_{C^2(M)} \le M \hbox{ and } \inf_M u \ge \frac{1}{M} \right \},
\end{equation*}
and let $\theta^{-} <  \theta^{+}$ be such that equation \eqref{ELt} has at least a solution for $\theta = \theta^{-}$ and has none for $\theta = \theta^{+}$. Define, for any $\theta^{-} \le \theta \le \theta^{+}$, $J_\theta: \{ u \in C^2(M), \inf_M u  > 0\} \to C^2(M)$ by 
\begin{equation*}
J_\theta(u) = u - \left( \triangle_g + h \right)^{-1} \left( f u^{2^*-1} + \frac{\theta a}{u^{2^*+1}} \right) .
\end{equation*}
By Theorem \ref{Thstabi}, there exists a positive $M_0$ such that $J_\theta^{-1}(\{ 0 \}) \subset \Omega_{M_0}$ for all $\theta^{-} \le \theta \le \theta^{+}$, so that the Leray-Schauder degree $\textrm{deg} (J_\theta, \Omega_{M}, 0)$ is well-defined for $\theta^{-} \le \theta \le \theta^{+}$ and for any $M > M_0$. Since it is homotopy-invariant and no solutions of \eqref{ELt} exist for $\theta = \theta^{+}$, there holds
\begin{equation} \label{calculdegre}
\textrm{deg} (J_\theta, \Omega_{M}, 0) = 0
\end{equation}
for any $M > M_0$ and any $ \theta^{-} \le \theta \le \theta^{+}$. In particular, \eqref{calculdegre} shows that solutions of \eqref{ELt} generically appear by pairs. Theorem \ref{Th1} shows that two solutions actually always exist. 

\medskip  Existence of at least one solution for \eqref{ELFull} without the CMC-assumption can be found in Allen-Clausen-Isenberg \cite{ACI}, 
Dahl-Gicquaud-Humbert \cite{DaGiHu}, Holst-Nagy-Tsogtgerel \cite{HoNaTso} and Maxwell \cite{Maxwell} when $\mathcal{B}_{\tau,\psi,V} < 0$, and in 
Premoselli \cite{Premoselli} when $\mathcal{B}_{\tau,\psi,V} > 0$. 

\section{Sharp blow-up estimates.} \label{stabilite}

We let $(M,g)$ be a smooth closed Riemannian manifold of dimension $3 \le n \le 5$, $h,f$ and $a$ be smooth functions on $M$ satisfying the assumptions of Theorem \ref{Thstabi}. 
We recall that $\triangle_g + h$ is coercive if there exists a positive constant $C$ such that for any $u \in H^1(M)$,
\[  \int_M \left( |\nabla u|_g^2 + hu^2 \right) dv_g \geqslant C \| u \|_{H^1(M)}^2 \]
or, equivalently, if
\begin{equation} \label{Nh}
\Vert u \Vert_{H^1_h} =  \left( \int_M \left( |\nabla u|_g^2 + hu^2 \right) dv_g \right)^{\frac{1}{2}}
\end{equation}
is an equivalent norm on $H^1(M)$. In this case we define the constant $S_h$ as the smallest positive constant satisfying for all $u \in H^1(M)$:
\begin{equation} \label{Sh}
\Vert u \Vert_{L^{2^*}} \le S_{h}^{\frac{1}{2^*}} \Vert u \Vert_{H^1_h}.
\end{equation}
We let $(a_k)_k$ be a sequence of non-negative functions converging to $a$ in $C^0(M)$ as $k \to \infty$ and 
$(q_k)_k$, $2 \le q_k \le 2^*$, be a sequence of positive real numbers converging to $2^*$. Let $(u_k)_k$ be a sequence of solutions of \eqref{ELq}. We assume that the following assumption holds:
\begin{equation} \label{grossehypo}
\textrm{ either } f > 0 \quad \textrm{ or } \quad  \| u_k \|_{H^1} \le C_0\hskip.1cm ,
\end{equation}
where $C_0$ is a positive constant independent of $k$. As a first striking fact, we prove that the sequence $(u_k)_k$ is uniformly bounded from below by some positive number:

\begin{prop}
Let $(u_k)_k$ be a sequence of solutions of \eqref{ELq}. There exists $\ve_0 > 0$ independent of $k$ such that
\begin{equation} \label{unifminor}
u_k \ge \ve_0
\end{equation}
for all $k$.
\end{prop}

\begin{proof}
We follow the arguments in Hebey-Veronelli \cite{HeVer}.  Given $K > 0$, we let $H = h + K $ and choose $K$ large enough 
in order to have $H \ge 1$. For any $\delta > 0$ we consider the unique functions $\psi_{k,\delta}$, $\psi_\delta$ and $\psi_0$ solving:
\begin{equation} \label{astuceminor}
\left \{ 
\begin{aligned}
& \triangle_g \psi_{k, \delta} + H \psi_{k, \delta} = a_k - \delta f^{-}\hskip.1cm , \\
& \triangle_g \psi_\delta + H \psi_{\delta} = a - \delta f^{-}\hskip.1cm , \\
& \triangle_g \psi_0 + H  \psi_{0} = a \hskip.1cm ,
\end{aligned} \right.
\end{equation}
where we have let $f^{-} = - \min(f,0)$. By standard elliptic theory, $\psi_{k,\delta} \to \psi_\delta$ in $C^0(M)$ as $k \to \infty$ and $\psi_\delta \to \psi_0$ in $C^0(M)$ as $\delta \to 0$.  Since $a \not \equiv 0$, by the maximum principle there holds $\psi_0 > 0$ in $M$ and so, for some $\delta_0$ small enough, $\psi_{\delta_0} >0$ in $M$. In particular, $\psi_{k, \delta_0} \ge \frac{1}{2} \psi_{\delta_0}$ for $k$ large enough. Consider now $\theta_k = t \psi_{k, \delta_0}$, $t >0$. Since $q_k \to 2^*$ and $a_k \to a$ in $C^0(M)$ as $k \to \infty$, for $t$ sufficiently small (that does not depend on $k$) there holds for any $k$ large enough:
\[ \triangle_g \theta_k + H \theta_k \le \frac{a_k}{\theta_k^{q_k+1}} - f^{-} \theta_k^{q_k-1}.  \]
Hence:
\begin{equation}
\begin{aligned}
  \triangle_g (u_k - \theta_k) +  H  (u_k - \theta_k) & \ge f u_k^{q_k-1} + f^{-} \theta_k^{q_k-1} + a_k \left( u_k^{-q_k-1} - \theta_k^{-q_k-1} \right) \\
  & \ge 0
\end{aligned}
\end{equation}
at any point $x\in M$ such that $u_k(x) \le \theta_k(x)$. By the maximum principle, we thus have $u_k \ge \theta_k$ in $M$. Since $\theta_k \ge \frac{1}{2} t \psi_{\delta_0}$, we have a uniform lower bound for the $u_k$.
\end{proof}

With \eqref{unifminor} we get in particular the existence of some constant $C$ depending on $h,f,a$ and $\ve_0$ such that
\begin{equation} \label{bornelapla}
| \triangle_g u_k| \le C u_k^{q_k-1}.
\end{equation}
Standard elliptic theory and \eqref{bornelapla} show that Theorem \ref{Thstabi} is proved provided $(u_k)_k$ is uniformly bounded in $L^\infty(M)$. We thus proceed by contradiction and assume that 
\begin{equation} \label{blowup}
\max_M u_k \to +\infty
\end{equation} 
as $k \to \infty$. 
In this section we perform an asymptotic analysis of the sequence $(u_k)_k$ around a concentration point and obtain sharp pointwise estimates.  We denote the injectivity radius of $(M,g)$ by $i_g$. Following Druet-Hebey \cite{DruHeb}, we let  $(x_k)_k$ be a sequence of points in $M$ and $(\rho_k)_k$ be a sequence of positive numbers, with $ 0 < \rho_k < i_g / 7$, satisfying for any $k \in \mathbb{N}$:
\begin{equation} \label{hyposuites}
\left \{
\begin{aligned}
& x_k \textrm{ is a critical point of } u_k\hskip.1cm , \\
& d_g(x_k,x)^{\puiq} u_k(x) \le C \textrm{ for all } x \in B_{x_k}(7 \rho_k)\hskip.1cm ,  \\
& \lim_{k \to \infty} \rho_k^{\puiq} \sup_{B_{x_k}(6 \rho_k)} u_k = + \infty .\\ 
\end{aligned}
\right.
\end{equation}
In \eqref{hyposuites} and in everything that follows we will denote by $C$ some constant that does not depend on $k$. We define $\mu_k$ as
\begin{equation} \label{muq}
\mu_k = u_k(x_k)^{- \frac{q_k-2}{2}}.
\end{equation}
We can apply here the asymptotic analysis of equation \eqref{ELq} and of subcritical perturbations of the Yamabe equation as found respectively in Druet-Hebey \cite{DruHeb} and Druet \cite{DruYlowdim}. Similar a priori blow-up techniques were first developed for the analysis of compactness of solutions of the Yamabe equation by Schoen \cite{SchoenPreprint}, Li-Zhu \cite{LiZhu} and by Druet \cite{DruYlowdim}. Assuming \eqref{hyposuites} and \eqref{grossehypo} we obtain that 
\[ \mu_k^{- \frac{2}{q_k-2}} \sim \sup_{B_{x_k}(6 \rho_k)} u_k, \quad  \frac{\rho_k}{\mu_k} \to + \infty   \]
(so, in particular, $\mu_k \to 0$ as $k \to \infty$), that
\begin{equation} \label{fpositif}
f(x_0) > 0\hskip.1cm ,
\end{equation}
where $x_0$ is a limit of a subsequence of $x_k$, and that
\begin{equation} \label{convC1locq}
\mu_k^{\puiq} u_k \big( \exp_{x_k} (\mu_k x ) \big) \to \left( 1 + \frac{f(x_0)}{n(n-2)} |x|^2 \right)^{- \pui} 
\end{equation}
in $C^1_{loc}(\RR^n)$ as $k\to \infty$, where $\mu_k$ is as in \eqref{muq}. Also the Harnack inequality stated in Druet-Hebey \cite{DruHeb} (lemma $1.3$) is still satisfied here: there exists $ C > 1 $ such that for any sequence $(s_k)_k$ of positive real numbers with $0 < 6 s_k \leq \rho_k$ we have
\begin{equation} \label{Harnackq}
s_k ||\nabla u_k ||_{L^\infty(\Omega_k)} \leq C \sup_{\Omega_k} u_k \leq C^2 \inf_{\Omega_k} u_k\hskip.1cm ,
\end{equation}
where $\Omega_k = B_{x_k}(6 s_k) \backslash B_{x_k}(\frac{1}{6}s_k)$. Now we define $\vp_k: (0,\rho_k) \mapsto \mathbb{R}^+$ as the mean value of $u_k$ on spheres centered at $x_k$:

\begin{equation} \label{defvpq}
 \vp_k(r) = \frac{1}{| \partial B_{x_k}(r)|_g} \int_{\partial B_{x_k}(r)} u_k d\sigma_g .
 \end{equation}
Also we define $r_k$, the radius of influence of the bubble centered at $x_k$:
\begin{equation} \label{proprq1}
r_k = \sup \left\{ r \in (2R_0 \mu_k ; \rho_k) \textrm{ such that } s^{\frac{2}{q_k-2}}\vp_k(s) \textrm{ is nonincreasing in } (2 R_0 \mu_k ; r) \right\}\hskip.1cm ,
\end{equation}
where we have let 
\begin{equation} \label{defR0}
R_0^2 = \frac{n(n-2)}{f(x_0)}.
\end{equation} 
With \eqref{convC1locq}, for any positive $R$ there holds, as $k \to \infty$:
\begin{equation} \label{borddebulle}
(R \mu_k)^{\frac{2}{q_k-2}} \vp_k(R \mu_k) \to R^{\pui} \left(1 + \frac{R^2}{R_0^2} \right)^{1 - \frac{n}{2}}\hskip.1cm ,
\end{equation}
where $R_0$ is as in \eqref{defR0}. Hence
\begin{equation} \label{rsurmu}
\frac{r_k}{\mu_k} \to + \infty.
\end{equation}
Note that from the definition of $r_k$ there holds that
\begin{equation} \label{proprq2}
r^{\frac{2}{q_k-2}} \vp_k \textrm{ is nonincreasing in } (2 R_0 \mu_k, r_k)
\end{equation}
and that
\begin{equation} \label{proprq3}
\textrm{ if } r_k < \rho_k, \left(r^{\frac{2}{q_k-2}} \vp_k(r) \right)'(r_k) = 0.
\end{equation}
As consequence of \eqref{Harnackq} one has that for $0 < s \le r_k$ and for some positive constant $C$:
\begin{equation} \label{lemme51}
\frac{1}{C} \sup_{B_{x_k}(6 s) \backslash B_{x_k}(\frac{1}{6} s)} u_k \le \vp_k(s) \le C \inf_{B_{x_k}(6 s) \backslash B_{x_k}(\frac{1}{6} s)} u_k. 
\end{equation}
In particular, by \eqref{proprq2}, there holds for any $R > 0$:
\[ \sup_{x \in B_{x_k}(6r_k) \backslash B_{x_k}(R \mu_k)} d_g(x_k,x)^{\frac{2}{q_k-2}} u_k(x) \le (R \mu_k)^{\frac{2}{q_k-2}} \vp_k(R \mu_k)\]
so that, with \eqref{borddebulle}:
\begin{equation} \label{estimeefaible}
\lim_{R \to \infty} \lim_{k \to \infty} \sup_{x \in B_{x_k}(6r_k) \backslash B_{x_k}(R \mu_k)} d_g(x_k,x)^{\frac{2}{q_k-2}} u_k(x) = 0.
\end{equation}
Let us define
\begin{equation} \label{etaq}
 \eta_k = \sup_{B_{x_k}(6 r_k) \backslash B_{x_k}(\frac{1}{6}r_k)} u_k .
\end{equation}
In particular \eqref{Harnackq}, \eqref{estimeefaible} and \eqref{unifminor} imply that:
\begin{equation} \label{taillerqetaq}
 r_k^2 \eta_k^{q_k-2} \to 0 \quad \textrm{and} \quad r_k \to 0 .
 \end{equation}
Now we prove a sharp asymptotic control from above on $u_k$ on a ball of radius $r_k$.

\begin{lemme} \label{lemma6}
Let $(u_k)_k$ be a sequence of solutions of \eqref{ELq} satisfying \eqref{grossehypo} and \eqref{blowup}. Let  $(x_k)_k$ and $(\rho_k)_k$ be such that \eqref{hyposuites} holds. There exists a positive constant $C$ such that, for any $ x \in B_{x_k}(6 r_k) \backslash \{x_k\}$:
\begin{equation} \label{optimalq}
 u_k(x) + d_g(x_k,x) | \nabla u_k (x) | \le C \mu_k^{n - \frac{2(q_k-1)}{q_k-2}} d_g(x_k,x)^{2-n} .
 \end{equation}
As a consequence:
\begin{equation} \label{taillerq}
 r_k^{n-2} = O \left( \mu_k^{n - \frac{2(q_k-1)}{q_k-2}} \right)\hskip.1cm ,
 \end{equation}
 where $\mu_k$ is as in \eqref{muq}.
\end{lemme}

\begin{proof} We first show that there exists a positive constant $C$ such that
\begin{equation} \label{optimalqeta}
u_k(x) \le C \left(  \mu_k^{n-\frac{2(q_k-1)}{q_k-2}} d_g(x_k,x)^{2-n} + \eta_k \right)
\end{equation}
for all $x \in B_{x_k}(6 r_k) \backslash \{x_k\}$ and for all $k$, where $\eta_k$ is as in \eqref{etaq}. To prove this we let $(y_k)_k$ be an arbitrary sequence in $B_{x_k}(6 r_k) \backslash \{x_k\}$ and prove that there exists a positive constant $C$ such that, up to a subsequence:
\begin{equation} \label{optimalqetasuite}
u_k(y_k) = O \left( \mu_k^{n - \frac{2(q_k-1)}{q_k-2}} d_g(x_k,y_k)^{2-n} \right) + O(\eta_k).
\end{equation}
First, \eqref{optimalqetasuite} follows from \eqref{convC1locq} if $d_g(x_k,y_k) = O(\mu_k)$ and from \eqref{Harnackq} if $\frac{d_g(x_k,y_k)}{r_k} \not \to 0$ when $k \to \infty$. We thus assume from now on that
\begin{equation} \label{lemme63}
\frac{d_g(x_k,y_k)}{\mu_k} \to + \infty \hbox{ and } \frac{d_g(x_k,y_k)}{r_k} \to 0 .
\end{equation}
We let $G_k$ be the Green function of $\triangle_g$ in the ball $B_{x_k}(6 r_k)$ with Dirichlet boundary condition satisfying $G \ge 1$. We recall that there exists a continuous function $\tau: \mathbb{R}^+ \to \mathbb{R}^+$, $\tau(0)=0$ such that
\begin{equation} \label{Green1}
\left| G_k(x,y) - \frac{1}{(n-2)\omega_{n-1}d_g(x,y)^{n-2}} \right| \le \tau\left(d_g(x,y)\right) d_g(x,y)^{2-n}
\end{equation}
and
\begin{equation} \label{Green2}
\left| |\nabla G_k(x,y)| - \frac{1}{\omega_{n-1}d_g(x,y)^{n-1}} \right| \le \tau\left(d_g(x,y)\right) d_g(x,y)^{1-n}
\end{equation}
(see Aubin \cite{Aub} or Robert \cite{Robert} for such estimates). The representation formula for the Green function gives, with \eqref{bornelapla}, \eqref{lemme63} and \eqref{Green1}:
\begin{equation} \label{representation}
\begin{aligned}
u_k(y_k) & = \int_{B_{x_k}(6 r_k)} G_k(x,y_k) \triangle_g u_k(x) dv_g(x) - \int_{\partial B_{x_k}(6 r_k)} \partial_{\nu} G_k(x,y_k) u_k(x) d\sigma_g(x) \\
& = O \left( \int_{B_{x_k}(6 r_k)} d_g(x,y_k)^{2-n} u_k^{q_k-1}(x) dv_g(x) \right) + O(\eta_k) .\\
\end{aligned} 
\end{equation}
Combining the subcritical analysis in Druet \cite{DruYlowdim} and the analysis for the critical equation in Druet-Hebey \cite{DruHeb} one gets that for any $\ve > 0$ there exists a positive constant $C_\ve$ such that
\begin{equation} \label{optimalveq}
 u_k(x) \le C_\varepsilon \left( \mu_k^{\frac{2}{q_k-2}(1-2\varepsilon)}d_g(x_k,x)^{-(1-\varepsilon)\frac{4}{q_k-2}} + \eta_k \left( \frac{r_k}{d_g(x_k,x)}\right)^{\frac{4}{q_k-2}\varepsilon} \right) 
 \end{equation}
for all $ x \in B_{x_k}(6 r_k) \backslash \{x_k\}$. Coming back to \eqref{representation}, we have
\begin{equation} \label{ints}
 \begin{aligned}
\int_{B_{x_k}(6 r_k)} d_g(x,y_k)^{2-n} u_k^{q_k-1}(x) dv_g(x) & = \int_{B_{x_k}(\mu_k)} d_g(x,y_k)^{2-n} u_k^{q_k-1}(x) dv_g(x)\\
 \quad & + \int_{B_{x_k}(6 r_k) \backslash B_{x_k}(\mu_k)} d_g(x,y_k)^{2-n} u_k^{q_k-1}(x) dv_g(x) .\\ 
\end{aligned} 
\end{equation}
One finds, on one hand using \eqref{convC1locq} and \eqref{lemme63} that
\begin{equation} \label{int1}
\int_{B_{x_k}(\mu_k)} d_g(y_k,x)^{2-n} u_k^{q_k-1}(x) dv_g(x) = O \left( \mu_k^{n - \frac{2(q_k-1)}{q_k-2}} d_g(x_k,y_k)^{2-n} \right)
\end{equation}
and on the other hand using \eqref{optimalveq} and \eqref{lemme63} that
\begin{equation} \label{int2}
 \begin{aligned}
& \int_{B_{x_k}(6 r_k) \backslash B_{x_k}(\mu_k)} d_g(y_k,x)^{2-n} u_k^{q_k-1}(x) dv_g(x) \\
& \quad  = O \left( \mu_k^{\frac{2(q_k-1)}{q_k-2}(1-2\varepsilon)} \int_{B_{x_k}(6 r_k) \backslash B_{x_k}(\mu_k)} d_g(y_k,x)^{2-n} d_g(x_k,x)^{-\frac{4(q_k-1)}{q_k-2}(1-\varepsilon)} dv_g(x) \right) \\
& \qquad + O \left( \int_{B_{x_k}(6 r_k) \backslash B_{x_k}(\mu_k)} \eta_k^{q_k-1} r_k^{\frac{4(q_k-1)}{q_k-2} \varepsilon} d_g(y_k,x)^{2-n} d_g(x_k,x)^{-\frac{4(q_k-1)}{q_k-2} \varepsilon} dv_g(x) \right) \\
& = O \left( \mu_k^{n - \frac{2(q_k-1)}{q_k-2}} d_g(x_k,y_k)^{2-n} \right) + O \left( r_k^2 \eta_k^{q_k-1} \right) \\
& = O \left( \mu_k^{n - \frac{2(q_k-1)}{q_k-2}} d_g(x_k,y_k)^{2-n} \right)  + o(\eta_k) \hskip.1cm ,
\end{aligned}
\end{equation}
where we used \eqref{taillerqetaq} to obtain the last equality in \eqref{int2}. In the end, \eqref{representation},  \eqref{ints} , \eqref{int1} and \eqref{int2} yield \eqref{optimalqetasuite} and in particular \eqref{optimalqeta} holds true. Now we show that there exists a positive constant $C$ such that for $k$ large enough: 
\begin{equation} \label{ineqetaq}
\eta_k \le C \mu_k^{n - \frac{2(q_k-1)}{q_k-2}} r_k^{2-n}.
\end{equation}
Using \eqref{proprq2} we can assert that for any $0 < \delta < 1$:
\begin{equation} \label{ineqdelta}
(\delta r_k)^{\frac{2}{q_k-2}} \vp_k(\delta r_k) \ge r_k^{\frac{2}{q_k-2}} \vp_k(r_k) .
\end{equation}
So \eqref{lemme51}, \eqref{ineqdelta}  and the definition of $\vp_k$ as in \eqref{defvpq} give, for some positive constant $C$:
\[ \frac{1}{C} r_k^{\frac{2}{q_k-2}} \eta_k \le r_k^{\frac{2}{q_k-2}} \vp_k(r_k) \le (\delta r_k)^{\frac{2}{q_k-2}} \sup_{\partial B_{x_k}(\delta r_k)} u_k\hskip.1cm , \]
where $\eta_k$ is as in \eqref{etaq}. Using \eqref{optimalqeta} one gets for some positive $C$:

\[ \frac{1}{C} \eta_k \le \delta^{\frac{2}{q_k-2}}  \left( \mu_k^{n - \frac{2(q_k-1)}{q_k-2}} (\delta r_k)^{2-n} + \eta_k \right) .\] 
We just need to choose $\delta$ small enough to have $ C \delta^{\frac{2}{q_k-2}} < 1 $ to obtain 
\[ \eta_k \le C \mu_k^{n - \frac{2(q_k-1)}{q_k-2}} (\delta r_k)^{2-n} \]
for some positive $C$, which shows \eqref{ineqetaq}. Then \eqref{ineqetaq} and  \eqref{optimalqeta} conclude the sharp asymptotic \eqref{optimalq} for $u_k$. The gradient part in \eqref{optimalq} is obtained choosing $s_k = d_g(x_k,y_k)$, for any sequence $y_k \in B_{x_k}(6 r_k)$, in \eqref{Harnackq}. Finally \eqref{taillerq} is obtained using \eqref{unifminor} in \eqref{optimalq}, and this concludes the proof of Lemma \ref{lemma6}
\end{proof}

With Lemma \ref{lemma6} we now determine the asymptotical behavior of the sequence $(u_k)_k$ at the boundary of its maximal ball of influence:

\begin{lemme} \label{lemma7}
Let $(u_k)_k$ be a sequence of solutions of \eqref{ELq} satisfying \eqref{grossehypo} and \eqref{blowup}. Let  $(x_k)_k$ and $(\rho_k)_k$ be such that \eqref{hyposuites} holds true. Then up to a subsequence, we have $\rho_k = r_k$, $\rho_k \to 0$, $\rho_k^{\frac{2}{q_k-2}} u_k(x_k) \to + \infty$ and 
there exists a harmonic function $H$ in $B_0(5)$ with $H(0) \le 0$ such that
\begin{equation} \label{asympto}
 \mu_k^{\frac{2(q_k-1)}{q_k-2}-n} \rho_k^{n-2} u_k\left( \exp_{x_k}(\rho_k x) \right) \to \frac{R_0^{n-2}}{|x|^{n-2}} + H(x) 
 \end{equation}
in $C^2_{loc}(B_0(5) \backslash \{0\})$ as $k \to \infty$, where $R_0$ is as in \eqref{defR0}.
\end{lemme}
\begin{rem} 
With the definition of $\mu_k$ in \eqref{muq} one can also write \eqref{asympto} as 
\begin{equation} \label{remasympto}
 u_k(x_k)^{1 -n + q_k \frac{n-2}{2}} \rho_k^{n-2} u_k \left( \exp_{x_k}(\rho_k x) \right) \to \frac{R_0^{n-2}}{|x|^{n-2}} + H(x). 
 \end{equation}
\end{rem}

\begin{proof}
For $x \in B_0(5)$ we define
\[ \begin{aligned}
 \hat{u}_k(x) & = \mu_k^{\frac{2(q_k-1)}{q_k-2}-n} r_k^{n-2} u_k \left( \exp_{x_k}(r_k x)\right), \\
  \hat{g}_k(x)&  = \left( \exp_{x_k}^* \right) g(r_k x). \\
 \end{aligned} 
\]
Since $r_k \to 0$ by \eqref{taillerqetaq}, $\hat{g}_k \to \xi$ the euclidean metric in $C^2_{loc}(B_0(5))$. Thanks to \eqref{optimalq} we also have that 
\begin{equation} \label{singhatuq}
\hat{u}_k(x) \le C |x|^{2-n} 
\end{equation}
in $B_0(5) \backslash \{ 0\}$. By \eqref{ELq} and \eqref{singhatuq} we therefore obtain that for all $k$,
\begin{equation} \label{lemme73}
\triangle_{\hat{g}_k} \hat{u}_k  = \hat{F}_k(x)
\end{equation}
in $B_0(5)  \backslash \{ 0\}$, where 
\begin{equation} \label{taillehatFq}
|\hat{F}_k(x)| \le  C \left( \frac{\mu_k}{r_k} \right)^{(n-2)(q_k-2)-2} |x|^{(q_k-1)(2-n)}
\end{equation}
for some positive $C$. By standard elliptic theory, \eqref{lemme73}, \eqref{taillehatFq} and by \eqref{rsurmu} there exists $\hat{U}$ harmonic in $B_0(5) \backslash \{0\}$ such that $\hat{u}_k \to \hat{U}$ in $C^2_{loc}(B_0(5) \backslash \{0\})$. Since by \eqref{singhatuq} $\hat{U}$ satisfies that $ 0 \le \hat{U}(x) \le C |x|^{2-n}$ for some positive constant $C$ in $B_0(5) \backslash \{0\}$ we can write
\begin{equation} \label{decompohatu}
 \hat{U}(x) = \frac{\lambda}{|x|^{n-2}} + H(x)\hskip.1cm ,
 \end{equation}
where $\lambda \ge 0$ and $H$ is harmonic in $B_0(5)$. To show that $\lambda = R_0^{n-2}$ we integrate \eqref{lemme73} on $B_0(1)$: there holds
\begin{equation} \label{lemme74}
 - \int_{\partial B_0(1)} \partial_{\nu} \hat{u}_k d \sigma_{\hat{g}_k} = \int_{B_0(1)} \hat{F}_k dv_{\hat{g}_k} .
\end{equation}
Straightforward calculation gives, using \eqref{convC1locq} and \eqref{optimalq}:
\begin{equation} \label{rhslambda}
\begin{aligned}
\int_{B_0(1)} \hat{F}_k dv_{\hat{g}_k}  = & \int_{B_0(1)} \left( \frac{\mu_k}{r_k} \right)^{(n-2)(q_k-2)-2} f \big( \exp_{x_k}(r_k x) \big) \hat{u}_k^{q_k-1} (x) dx + o(1) \\
& = f(x_0) \int_{\mathbb{R}^n} \left( 1+ \frac{|x|^2}{R_0^2} \right)^{-1-\frac{n}{2}} dx + o(1).\\
& = (n-2) \omega_{n-1} R_0^{n-2}  + o(1) .\\
\end{aligned}
\end{equation}
On the other side \eqref{decompohatu} shows that  
\begin{equation} \label{lhslambda}
\int_{\partial B_0(1)} \partial_{\nu} \hat{u}_k d\sigma_{\hat{g}_k} \to - \lambda (n-2) \omega_{n-1} 
\end{equation}
as $k\to\infty$. So \eqref{rhslambda} and \eqref{lhslambda} in \eqref{lemme74} yield $\lambda = R_0^{n-2}$. In the end,
\begin{equation} \label{convdistrq}
\mu_k^{\frac{2(q_k-1)}{q_k-2}-n} r_k^{n-2} u_k \left( \exp_{x_k}(r_k x)\right) \to  \hat{U} = \frac{R_0^{n-2}}{|x|^{n-2}} + H(x) 
\end{equation}
in $C^2_{loc}(B_0(5) \backslash \{0 \})$ as $k \to \infty$. We now prove that $r_k = \rho_k$. Assume first by contradiction that $r_k < \rho_k$. Let for any $r > 0$:
\[ \varphi(r) = \frac{1}{\omega_{n-1}r^{n-1}} \int_{\partial B_0(r)} \hat{U} d \sigma\hskip.1cm , \]
where $\hat{U}$ is as in \eqref{decompohatu}. On one side, \eqref{convdistrq} shows that 
\begin{equation} \label{valvp}
\vp(r) = \left( \frac{R_0}{r} \right)^{n-2} + H(0). 
\end{equation}
On the other side \eqref{proprq3} and \eqref{convdistrq} give us that
\begin{equation} \label{vpdiff1}
 \left( r^{\frac{n-2}{2}} \varphi(r) \right)'(1) = 0 .
 \end{equation}
 Combining \eqref{valvp} and \eqref{vpdiff1} yields 
\[ H(0) = R_0^{n-2} .\]
To conclude the proof of Lemma \ref{lemma7} and show that $r_k = \rho_k$ we therefore show that $H(0) \le 0$.  We let $X_k(x) = \nabla f_k(x)$, where $f_k(x) = \frac{1}{2} d_g(x_k,x)^2$, be the vector field whose coordinates in the exponential chart at $x_k$ are $x^i$ and use the Pohozaev identity as in Druet-Hebey \cite{DruHeb} (see also Druet \cite{DruYlowdim}) on $\Omega_k = B_{x_k}(6 r_k)$: 
\begin{equation} \label{Pohozaev}
\begin{aligned}
& \int_{\Omega_k} \nabla u_k (X_k) \triangle_g u_k dv_g  = \int_{\Omega_k} \left( \nabla^{\sharp} X_k ( \nabla u_k, \nabla u_k) - \frac{1}{2} (\textrm{div}_g X_k) |\nabla u_k|^2 \right) dv_g \\
& + \int_{\partial \Omega_k} \left( \frac{1}{2}(X_k, \nu)_g |\nabla u_k|^2 - \nabla u_k(X_k) \partial_{\nu} u_k \right) d \sigma_g\hskip.1cm , \\
\end{aligned}
\end{equation} 
where we have let $\left(\nabla^{\sharp} X_k\right)^{ij} = \left( \nabla^i X_k \right)^{j}$, where $(\cdot, \cdot)_g$ is the scalar product defined by $g$ and $\nu$ is the unit outward normal to $\partial \Omega_k$. 
In the exponential chart at $x_k$ we have 
\begin{equation} \label{carteexpo}
\left( \nabla^i X_k \right)^j = g^{ij} + O \left( d_g(x_k,x)^2 \right) 
\end{equation}
and so integrating by parts the first integral in the right-hand side of \eqref{Pohozaev} and using \eqref{carteexpo} gives
\begin{equation} \label{egalpoho}
\begin{aligned}
& \int_{\Omega_k} \left( \nabla u_k(X_k) + \frac{n-2}{2} u_k \right) \triangle_g u_k dv_g =  O \left( \int_{\Omega_k} d_g(x_k,x)^2 |\nabla u_k|^2 dv_g \right) \\
& + \int_{\partial \Omega_k} \left( \frac{1}{2} (X_k, \nu)_g |\nabla u_k|^2 - \nabla u_k (X_k) \partial_{\nu} u_k - \frac{n-2}{2} u_k \partial_{\nu} u_k \right) d \sigma_g
\hskip.1cm .\\
\end{aligned}
\end{equation}
Using \eqref{convdistrq} it is easily seen that, since $H$ is harmonic:
\begin{equation} \label{intomegaqstep}
\begin{aligned}
& \mu_k^{\frac{4(q_k-1)}{q_k-2}-2n} r_k^{n-2} \int_{\partial \Omega_k} \left( \frac{1}{2} (X_k, \nu)_g |\nabla u_k|^2 - \nabla u_k (X_k) \partial_{\nu} u_k - \frac{n-2}{2} u_k \partial_{\nu} u_k \right) d \sigma_g \\
 & =  \int_{\partial B_0(1)} \left( \frac{1}{2} |\nabla \hat{U}|^2 - |\partial_{\nu} \hat{U}|^2 - \frac{n-2}{2} \hat{U} \partial_{\nu} \hat{U} \right) dx + o(1) \\
 & = \frac{(n-2)^2}{2} R_0^{n-2} \omega_{n-1} H(0) + \int_{\partial B_0(1)} \left( \frac{1}{2} |\nabla H|^2 - |\partial_\nu H|^2 - \frac{n-2}{2} H \partial_\nu H \right) dx + o(1) .
\end{aligned}
\end{equation}
A Pohozaev identity like \eqref{egalpoho} for $H$ and $X = x$ applied to the euclidean ball $B_0(1)$ for the euclidean metric gives, since $H$ is harmonic:
\begin{equation} \label{pohoH}
\begin{aligned}
 \int_{\partial B_0(1)} & \left( \frac{1}{2} |\nabla H|^2 - |\partial_\nu H|^2 - \frac{n-2}{2} H \partial_\nu H \right) dx \\
 & = \int_{B_0(1)} \left( \nabla H(X) + \frac{n-2}{2} H \right) \triangle_{\xi} H dx = 0.
 \end{aligned}
\end{equation}
We finally obtain, with \eqref{intomegaqstep} and \eqref{pohoH}:
\begin{equation} \label{intomegaq} 
\begin{aligned}
& \int_{\partial \Omega_k} \left( \frac{1}{2} (X_k, \nu)_g |\nabla u_k|^2 - \nabla u_k (X_k) \partial_{\nu} u_k - \frac{n-2}{2} u_k \partial_{\nu} u_k \right) d \sigma_g \\
& = \left( \frac{(n-2)^2}{2} R_0^{n-2} \omega_{n-1} H(0) + o(1) \right) \mu_k^{2n - \frac{4(q_k-1)}{q_k-2}} r_k^{2-n} .\\
\end{aligned} 
\end{equation}
We write now:
\[ \int_{\Omega_k} d_g(x_k,x)^2 |\nabla u_k|^2 dv_g = \int_{B_{x_k}(\mu_k)} d_g(x_k,x)^2 |\nabla u_k|^2 dv_g + \int_{B_{x_k}(r_k) \backslash B_{x_k}(\mu_k)} d_g(x_k,x)^2 |\nabla u_k|^2 dv_g \]
so that, using \eqref{convC1locq} and \eqref{optimalq} we obtain:
\begin{equation} \label{intdgnabla}
 \int_{\Omega_k} d_g(x_k,x)^2 |\nabla u_k|^2 dv_g = 
\left\{ 
 \begin{aligned}
 & O \left( \mu_k^{6 - \frac{4(q_k-1)}{q_k-2}} r_k \right)  & \textrm{ if } n = 3\hskip.1cm , \\
 & O \left( \mu_k^{8 - \frac{4(q_k-1)}{q_k-2}} \ln \frac{r_k}{\mu_k} \right) & \textrm{ if } n=4\hskip.1cm ,\\
 & O \left( \mu_k^{4+n - \frac{4(q_k-1)}{q_k-2}} \right) & \textrm{ if } n \ge 5\hskip.1cm . \\ 
 \end{aligned}
 \right. 
 \end{equation}
With \eqref{taillerq} equation \eqref{intdgnabla} becomes:
\begin{equation} \label{intomegaq2}
 \int_{\Omega_k} d_g(x_k,x)^2 |\nabla u_k|^2 dv_g = o \left( \mu_k^{2n- \frac{4(q_k-1)}{q_k-2}} r_k^{2-n} \right) 
 \end{equation}
for $3 \le n \le 5$. Now we compute the left-hand side integral in \eqref{egalpoho}. Using \eqref{ELq} and integrating by parts yield:

\begin{equation} \label{pohoavecf}
 \begin{aligned}
& \int_{\Omega_k} \left( \nabla u_k(X_k) + \frac{n-2}{2} u_k \right) \triangle_g u_k dv_g  
 = n \left( \frac{1}{2^*} - \frac{1}{q_k} \right) \int_{\Omega_k} f  u_k^{q_k} dv_g \\
 &  - \frac{1}{q_k} \int_{\Omega_k} \nabla f (X_k) u_k^{q_k} dv_g  + O \left( r_k \int_{\partial \Omega_k} u_k^{q_k} d \sigma_g \right)  + O \left( \int_{\Omega_k} d_g(x_k,x)^2 u_k^{q_k} dv_g \right) \\ 
 &  + O \left( \int_{\Omega_k} \big( d_g(x_k,x) |\nabla u_k| + u_k \big) \left( u_k^{-q_k-1} + u_k  \right) dv_g \right) .\\
\end{aligned} 
\end{equation}
Using \eqref{unifminor}, \eqref{convC1locq} and \eqref{optimalq}, the same computation as for \eqref{intdgnabla} gives, since $3 \le n \le 5$:
\begin{equation} \label{intpoho1} 
 \int_{\Omega_k} \big( d_g(x_k,x) |\nabla u_k| + u_k \big) \left( u_k^{-q_k-1} + u_k + 1\right) dv_g =   o\left( \mu_k^{2n - \frac{4(q_k-1)}{q_k-2}} r_k^{2-n} \right) .
 \end{equation}
Using \eqref{convdistrq} we obtain:
\begin{equation}  \label{intpoho2}
r_k \int_{\partial \Omega_k} u_k^{q_k} d \sigma_g = O \left( \left( \frac{\mu_k}{r_k}\right)^{(n-2)q_k} \mu_k^{- \frac{2 q_k}{q_k-2}} r_k^{n} \right) = o \left( \mu_k^{2n - \frac{4(q_k-1)}{q_k-2}} r_k^{2-n} \right). 
\end{equation}
Finally, using once again \eqref{convC1locq}, \eqref{optimalq} and \eqref{taillerq} there holds:
\begin{equation} \label{intpoho3}
\begin{aligned}
\int_{\Omega_k} d_g(x_k,x)^2 u_k^{q_k} dv_g & = O \left( \mu_k^{n- \frac{4}{q_k-2}} \right)  = o \left( \mu_k^{2n - \frac{4(q_k-1)}{q_k-2}} r_k^{2-n}\right)\hskip.1cm ,\\
\end{aligned}
\end{equation}
where the last equality is once again true because $3 \le n \le 5$. By \eqref{egalpoho} and \eqref{pohoavecf} and with \eqref{intomegaq}, \eqref{intomegaq2}, \eqref{intpoho1}, \eqref{intpoho2} and \eqref{intpoho3} we obtain:
\begin{equation} \label{derniereetape}
 \begin{aligned}
& \frac{1}{q_k} \int_{\Omega_k} \nabla f (X_k) u_k^{q_k} dv_g + n \left( \frac{1}{q_k} - \frac{1}{2^*} \right) \int_{\Omega_k} f  u_k^{q_k} dv_g \\
 & \quad  = \left( -\frac{(n-2)^2}{2} \omega_{n-1} R_0^{n-2} H(0) + o(1) \right) \mu_k^{2n - \frac{4(q_k-1)}{q_k-2}} r_k^{2-n}.
\end{aligned} \end{equation}
By \eqref{fpositif} we have $f(x) >0$ for $x \in \Omega_k$ for $k$ large enough so, since $q_k \le  2^*$ and $u_k \ge 0$, with \eqref{derniereetape} it is enough to show that
\begin{equation} \label{taillefpoho}
 \frac{1}{q_k} \int_{\Omega_k} \nabla f (X_k) u_k^{q_k} dv_g = o \left( \mu_k^{2n - \frac{4(q_k-1)}{q_k-2}}r_k^{2-n} \right) 
 \end{equation}
to conclude that $H(0) \le 0$. We first assume that $n=3$. For $R > 1 $ we write
\begin{equation} \label{sumintf}
 \int_{\Omega_k} \nabla f (X_k) u_k^{q_k} dv_g  = \int_{B_{x_k}(R \mu_k)} \nabla f (X_k) u_k^{q_k} dv_g + \int_{B_{x_k}(r_k) \backslash B_{x_k}(R \mu_k)} \nabla f (X_k) u_k^{q_k} dv_g  .
 \end{equation}
On one side, using \eqref{convC1locq} and since $x \mapsto x^i$ is an odd function in $\RR^n$ for any $i$, we have:
\begin{equation} \label{sumintf1}
   \int_{B_0(R \mu_k)} x^i \nabla_i f (x) u_k^{q_k}(x) dv_g(x)  = o \left( \mu_k^{n-1 - \frac{4}{q_k-2}}\right) .
\end{equation}
On the other side, using \eqref{optimalq} yields:
\begin{equation} \label{petitoR}
\begin{aligned} 
\int_{B_{x_k}(r_k) \backslash B_{x_k}(R \mu_k)} \nabla f (X_k) u_k^{q_k} dv_g & \le  C \mu_k^{n-1 - \frac{4}{q_k-2}} \int_{B_0(\frac{r_k}{\mu_k}) \backslash B_0(R)} |x|^{1+(2-n)q_k} dx  \\
& \le C \mu_k^{n-1 - \frac{4}{q_k-2} } R^{1+n +(2-n) q_k} \\
\end{aligned} 
\end{equation}
for some positive constant $C$ that depends neither on $R$ nor on $k$. Since $q_k \to 2^*$ as $k \to \infty$ \eqref{sumintf1} and \eqref{petitoR} show that: 
\begin{equation}  \label{sumintf2}
\begin{aligned}
\int_{\Omega_k} \nabla f (X_k) u_k^{q_k} dv_g & = o \left( \mu_k^{n-1 - \frac{4}{q_k-2} } \right)  = o \left( \mu_k^{6 - \frac{4(q_k-1)}{q_k-2}} r_k^{-1} \right)
\hskip.1cm ,
\end{aligned} 
\end{equation}
the last equality being true only when $n = 3$. So there holds $H(0) \le 0$ if $n=3$. 
Now we prove \eqref{taillefpoho} when $n=4,5$. Using \eqref{intpoho3} and \eqref{sumintf1} we have:
\begin{equation} \label{egalholder}
\begin{aligned}
\int_{\Omega_k} \nabla f (X_k) u_k^{q_k} dv_g & =  \partial_i f (x_k) \int_{\Omega_k} x^i u_k^{q_k} dv_g  + O \left( \int_{\Omega_k} d_g(x_k,x)^{2} u_k^{q_k} dv_g \right) \\
& = o \left( \mu_k^{n-1 - \frac{4}{q_k-2}} |\nabla f(x_k)| \right) +  o \left( \mu_k^{2n - \frac{4(q_k-1)}{q_k-2}} r_k^{2-n} \right) .\\
\end{aligned} 
\end{equation}
We now show that
\begin{equation} \label{taillenablafq}
 \mu_k^{n-1 - \frac{4}{q_k-2}} | \nabla f(x_k) | = 
o \left( \mu_k^{2n - \frac{4(q_k-1)}{q_k-2}} r_k^{2-n} \right)
\end{equation}
when $4 \le n \le 5$. With \eqref{egalholder} and \eqref{derniereetape} this will show that $H(0) \le 0$ when $4 \le n \le5$. To do this we apply another Pohozaev identity to $u_k$ in the open set $\Omega_k$. Let $Y \in \mathbb{R}^n$ be a given vector and let $Y_k$ be the local vector field whose coordinates in the exponential chart are $Y_k^i = Y^i$. Equality \eqref{Pohozaev} becomes:
\begin{equation} \label{lemme79}
\begin{aligned}
& \int_{\Omega_k} \nabla u_k (Y_k) \triangle_g u_k dv_g  = \int_{\Omega_k} \left( \nabla^{\sharp} Y_k (\nabla u_k, \nabla u_k) - \frac{1}{2} \textrm{div}_g Y_k |\nabla u_k|^2 \right) dv_g \\
& \quad + \int_{\partial \Omega_k} \left( \frac{1}{2} (Y_k, \nu)_g |\nabla u_k|^2 - \nabla u_k (Y_k) \partial_{\nu} u_k \right) d \sigma_g . \\
\end{aligned}
\end{equation}
Since $ \left( \nabla^i Y_k \right)^{j} = O \left( d_g(x_k,x) \right)$, \eqref{lemme79} becomes with \eqref{ELq}: 
\begin{equation} \label{lemme79bis}
 \begin{aligned}
& \int_{\Omega_k} \nabla u_k (Y_k) f u_k^{q_k-1} dv_g = O \left( \int_{\Omega_k} d_g(x_k,x) |\nabla u_k|^2 dv_g \right) \\
& \quad + O \left( \int_{\partial \Omega_k} |\nabla u_k|^2 d \sigma_g \right) + O \left( \int_{\Omega_k} |\nabla u_k| \left( u_k^{-q_k-1} + u_k  \right) dv_g \right) .\\
\end{aligned} 
\end{equation}
Using \eqref{convdistrq} we have
\[ \int_{\partial \Omega_k} |\nabla u_k|^2 d \sigma_g = O \left( \mu_k^{2n - \frac{4(q_k-1)}{q_k-2}} r_k^{1-n} \right) = o \left( \mu_k^{2n-1 - \frac{4(q_k-1)}{q_k-2}} r_k^{2-n} \right) .\]
Using \eqref{convC1locq}, \eqref{optimalq} and \eqref{taillerq} we find:
\[ \int_{\Omega_k} d_g(x_k,x) |\nabla u_k |^2 dv_g = O \left( \mu_k^{n-1 - \frac{4}{q_k-2}} \right) = o \left( \mu_k^{2n-1 - \frac{4(q_k-1)}{q_k-2}} r_k^{2-n} \right)\hskip.1cm , \]
the last equality being true since $n \le 5$. Once again, \eqref{unifminor}, \eqref{convC1locq}, \eqref{optimalq} and \eqref{taillerq} yield
\[ \int_{\Omega_k} |\nabla u_k| \left( u_k^{-q_k-1} + u_k + 1 \right) dv_g = O \left( \mu_k^{n-1- \frac{4}{q_k-2}} \right) =  o \left( \mu_k^{2n-1 - \frac{4(q_k-1)}{q_k-2}} r_k^{2-n}  \right) \]
since $4 \le n \le 5$. These computations with \eqref{lemme79bis} give: 
\begin{equation} \label{lemme710}
 \int_{\Omega_k} \nabla u_k (Y_k) f u_k^{q_k-1} dv_g = o \left( \mu_k^{2n-1 - \frac{4(q_k-1)}{q_k-2}} r_k^{2-n} \right) . 
\end{equation}
We now compute the latter integral differently: 
\[ \begin{aligned}
\int_{\Omega_k} \nabla u_k (Y_k) f u_k^{q_k-1} dv_g & = \frac{1}{q_k} \int_{\Omega_k} f \nabla ( u_k^{q_k} )(Y_k) dv_g \\
& = O \left( \int_{\partial \Omega_k} u_k^{q_k} d \sigma_g \right) - \frac{1}{q_k} \int_{\Omega_k} (\textrm{div}_g Y_k) f u_k^{q_k} dv_g \\
& \quad - \frac{1}{q_k} \int_{\Omega_k} \nabla f (Y_k) u_k^{q_k} dv_g \\
& = O \left( \int_{\partial \Omega_k} u_k^{q_k} d \sigma_g \right) + O \left( \int_{\Omega_k} d_g(x_k,x) u_k^{q_k} dv_g \right) \\
& \quad - \frac{1}{q_k} \left( \nabla f (Y_k) \right)(x_k) \int_{\Omega_k} u_k^{q_k} dv_g. \\
\end{aligned} \]
Using once again \eqref{convC1locq} and \eqref{optimalq} yields
\[ \int_{\Omega_k} d_g(x_k,x) u_k^{q_k} dv_g = O \left( \mu_k^{n-1 - \frac{4}{q_k-2}}\right) = o \left( \mu_k^{2n-1 - \frac{4(q_k-1)}{q_k-2}} r_k^{2-n} \right) \]
since $4 \le n \le 5$ and \eqref{intpoho2} gives
\[ \int_{\partial \Omega_k} u_k^{q_k} dv_g = o \left( \mu_k^{2n - \frac{4(q_k-1)}{q_k-2}} r_k^{1-n} \right)  = o \left( \mu_k^{2n-1 - \frac{4(q_k-1)}{q_k-2}} r_k^{2-n} \right). \]
Gathering the previous computation in \eqref{lemme710} gives in the end:
\begin{equation} \label{calculnablafq}
\frac{1}{q_k} \left( \nabla f (Y_k) \right)(x_k) \int_{\Omega_k} u_k^{q_k} dv_g = o \left( \mu_k^{2n-1 - \frac{4(q_k-1)}{q_k-2}} r_k^{2-n} \right) .
 \end{equation}
Finally, \eqref{convC1locq} and \eqref{optimalq} yield, with Lebesgue's dominated convergence theorem:
\[  \mu_k^{-n+2 + \frac{4}{q_k-2}} \int_{\Omega_k} u_k^{q_k} dv_g  =  \int_{\mathbb{R}^n} \left( 1 + \frac{|x|^2}{R_0^2} \right)^{-n} dx + o(1)\hskip.1cm ,  \]
so that with \eqref{calculnablafq} we obtain:
\begin{equation} \label{lafin}
  \mu_k^{n-2 - \frac{4}{q_k-2}} \left( \nabla f (Y_k) \right)(x_k) = o \left( \mu_k^{2n-1 - \frac{4(q_k-1)}{q_k-2}} r_k^{2-n} \right). 
  \end{equation}
Since  $Y_k^i = Y^i$ for all $i$, where $Y $ is an arbitrary vector in $\RR^n$, \eqref{lafin} implies \eqref{taillenablafq} and shows that $H(0) \le 0$ when $4 \le n \le 5$. In particular, this shows that $r_k = \rho_k$ and we have thus proven \eqref{asympto}. Noticing that then 
$$\rho_k^{\frac{2}{q_k-2}} u_k(x_k) = \left( \frac{r_k}{\mu_k} \right)^{\frac{2}{q_k-2}} \to + \infty$$
as $q_k \to 2^*$, this ends the proof of Lemma \ref{lemma7}.
\end{proof}

Note, as mentioned several times in the proof, that the above asymptotic analysis works only when $n \le 5$. 
When $n \ge 6$, counterexamples to Theorem \ref{Thstabi} for critical perturbations ($q_k = 2^*$ for all $k$) are known, see Druet-Hebey \cite{DruHeb}. 

\section{Stability under subcritical perturbations: proof of Theorem \ref{Thstabi}.} \label{preuvethstabi}

Using the asymptotic description obtained in the previous section we show that concentration points cannot appear and prove Theorem \ref{Thstabi}. Following the analysis in Druet-Hebey \cite{DruHeb}, see also Druet \cite{DruYlowdim}, there exists a positive constant $C$ such that for any $k $ there exists $N_k \in \mathbb{N}^*$ and $N_k$ critical points of $u_k$ denoted by $x_{1,k}, \dots, x_{N_k,k}$ such that 
\begin{equation} \label{ptscrit1}
d_g(x_{i,k}, x_{j,k})^{\frac{2}{q_k-2}} u_k(x_{i,k}) \geq 1 
\end{equation}
for all $i,j \in \{ 1, \dots, N_k\}$, $i\neq j$, and
\begin{equation} \label{ptscrit2}
 \left(\min_{i=1, \dots, N_k} d_g(x_{i,k},x) \right)^{\frac{2}{q_k-2}} u_k(x) \leq C_1 
 \end{equation}
for all $x \in M$ and for all $k$. Note that by construction, for any $k$ one finds among the $x_{i,k}$ the maximum point of $u_k$. As a consequence, $N_k \ge 2$ for $k$ large enough, otherwise thanks to \eqref{blowup} the hypothesis \eqref{hyposuites} would be satisfied with $\rho_k = \frac{1}{8} i_g(M)$ and with $x_k$ as the maximum point of $u_k$ which is impossible by Lemma \ref{lemma7}. Also note, with \eqref{ptscrit2}, that hypothesis \eqref{hyposuites} are satisfied if we choose $x_k = x_{i_k,k}$ for some $1 \le i_k \le N_k$ and for $\rho_k$ such that
\[ 7 \rho_k \le \min \left( \min_{1 \le i \le N_k, i \not = i_k} d_g(x_k, x_{i,k}) , \frac{1}{2}i_g(M) \right).\]
Let
\begin{equation} \label{defdq}
d_k = \min_{1 \le i < j \le N_k} d_g( x_{i,k}, x_{j,k}) 
\end{equation}
and assume, up to reordering the $x_{i,k}$, that 
\begin{equation} \label{propdq}
d_k = d_g(x_{1,k}, x_{2,k}) \le d_g(x_{1,k}, x_{3,k}) \le \cdots \le d_g(x_{1,k}, x_{N_k,k}).
\end{equation}
The following result shows that the concentration points are not isolated:

\begin{lemme} \label{lemma8}
Let $d_k$ be as in \eqref{defdq}. Then $d_k \to 0$ as $k \to \infty$.
\end{lemme}

\begin{proof}
We proceed by contradiction and assume that there exists $0 < d < i_g(M) $ such that for all $1 \le i < j \le N_k $ and for all $k$:
\begin{equation} \label{dtend0}
 d_g(x_{i,k},x_{j,k}) \ge d .
 \end{equation}
Since $M$ is compact the sequence $(N_k)_k$ is bounded and up to a subsequence we can assume that $N_k = N \ge 2$ for $k$ large enough. We let $x_k$ be a maximum point of $u_k$. By \eqref{blowup} and \eqref{ptscrit2} there exists $ i \in \{1, \cdots N \}$ such that $d_g(x_{i,k},x_k) \to 0 $ when $k \to \infty$. Hence \eqref{dtend0} shows that \eqref{hyposuites} is satisfied with $x_k$ and $\rho_k = \frac{1}{16}d$ which is impossible by Lemma \ref{lemma7}.
\end{proof}

Among the critical points $(x_{i,k})_{1 \le i \le N_k}$ we isolate those which are at a finite distance from $x_{1,k}$ when rescaling by a $d_k$ factor. For all $R>0$, we let $1 \le N_{R,k} \le N_k$ be such that
\begin{equation} \label{defNrk}
 \begin{aligned}
& d_g(x_{1,k},x_{i,k}) \le R d_k \textrm{ for } 1 \le i \le N_{R,k} \textrm{ and } \\
& d_g(x_{1,k},x_{i,k}) > R d_k  \textrm{ for } N_{R,k}+1 \le i \le N_k .\\
\end{aligned} 
\end{equation}
By \eqref{propdq} $N_{R,k}$ is well defined and we have $N_{R,k} = 1$ if $R<1$ and $N_{R,k} \ge 2$ otherwise. Let $R>1$. 
By \eqref{defdq} the balls $B_{x_{i,k}}(\frac{d_k}{4})$ and $B_{x_{j,k}}(\frac{d_k}{4})$ are disjoint and contained in $B_{x_{1,k}}((R+1)d_k)$ for all $1 \le i<j \le N_{R,k}$. Writing these inclusions in terms of volumes shows that, for any $R > 0$, $N_{R,k}$ is bounded by some positive constant that does not depend on $k$. To investigate the behavior of $u_k$ around $x_{1,k}$ let $ 0 < \delta \le \frac{1}{2} i_g(M)$ and define, for all $x \in B_0(\frac{\delta}{d_k})$:
\begin{equation} \label{checkuq}
\begin{split}
&\check{u}_k (x) = d_k^{\frac{2}{q_k-2}} u_k \left( \exp_{x_{1,k}}(d_k x) \right)\hskip.1cm ,\\
&\check{g}_k(x) = \left( \exp_{x_{1,k}}^* g \right)(d_k x)\hskip.1cm .
\end{split}
 \end{equation}
By Lemma \ref{lemma8}, $\check{g}_k \to \xi$ in $C^2_{loc}(\mathbb{R}^n)$. Using \eqref{ELq} we have in $B_0(\frac{\delta}{d_k})$:

\begin{equation} \label{eqcheckuq}
\triangle_{\check{g}_k} \check{u}_k = \check{F}_k\hskip.1cm ,
\end{equation}
where $\check{F}_k$ satisfies, by \eqref{bornelapla}:

\begin{equation} \label{propcheckfq}
| \check{F}_k | \le C_0 \check{u}_k^{q_k-1} .
\end{equation}
We finally let, for all $1 \le i \le N_k$ such that $d_g(x_{1,k}, x_{i,k}) \le \frac{1}{2} i_g(M)$:
\begin{equation} \label{checkxi}
 \check{x}_{i,k} = \frac{1}{d_k} \exp_{x_{1,k}}^{-1}(x_{i,k}) .
 \end{equation}
Let $R > 1$, $0 <r < 1$ and let
\begin{equation} \label{OmegarR}
 \Omega_{r,R} = B_0(R) \backslash \cup_{i=1}^{N_{2R,k}} B_{\check{x}_{i,k}}(r).
\end{equation}
Then \eqref{ptscrit2} shows that $\check{u}_k$ is uniformly bounded in $\Omega_{r,R}$ and mimicking the analysis of \eqref{Harnackq} we obtain a Harnack inequality for $\check{u}_k$ in $\Omega_{r,R}$: 
\begin{equation} \label{petitHarnack}
\Vert \nabla \check{u}_k \Vert_{L^\infty(\Omega_{r,R})} \le D_{r,R} \sup_{\Omega_{r,R}} \check{u}_k \le D_{r,R}^2 \inf_{\Omega_{r,R}} \check{u}_k\hskip.1cm ,
\end{equation}
where $D_{r,R} > 1$ is a constant depending only on $r$ and $R$. We have the following result that determines the behavior of $\check{u}_k$ around the critical points $\check{x}_{i,k}$ which are $d_k$-close to $\check{x}_{1,k}$.

\begin{lemme} \label{lemma10}
Let $(u_k)_k$ be a sequence of solutions of \eqref{ELq}, hence satisfying \eqref{unifminor}, such that \eqref{grossehypo} and \eqref{blowup} hold. Let $d_k$ be as in \eqref{defdq}, $\check{x}_{1,k}, \dots, \check{x}_{N_k,k}$ as in \eqref{checkxi} and $\check{u}_k$ as in \eqref{checkuq}. Then only one of the two following situations occurs: 
\begin{itemize}
\item either for any $1 \le i \le N_k$ such that $d_g(x_{1,k},x_{i,k}) = O(d_k)$ there holds
\begin{equation} \label{uqborne}
\sup_{B_{\check{x}_{i,k}}(\frac{1}{2})} \check{u}_k = O(1)
\end{equation}
as $k\to\infty$,  or
\item for any $1 \le i \le N_k$ such that $d_g(x_{1,k},x_{i,k}) = O(d_k)$
\begin{equation} \label{uqnonborne}
\sup_{B_{\check{x}_{i,k}}(\frac{1}{2})} \check{u}_k \to + \infty
\end{equation}
as $k \to \infty$.
\end{itemize}
\end{lemme}

\begin{proof}
Choose $k$ large enough and let $1 \le i \le N_k$ be such that $d_g(x_{1,k},x_{i,k}) = O(d_k)$. We denote by $\check{x}_i$ the limit of $(\check{x}_{i,k})_k$ in $\mathbb{R}^n$ when $k \to \infty$. Note that $|\check{x}_i| > 0$ if $i \ge 2$. We first investigate separately each of the cases stated in Lemma \ref{lemma10}. 
If \eqref{uqborne} holds then 
by \eqref{eqcheckuq}, \eqref{propcheckfq} and standard elliptic theory, $(\check{u}_k)_k$ is uniformly bounded in $C^1\left( B_{\check{x}_{i,k}}(\frac{1}{4}) \right)$. By \eqref{ptscrit1} we can thus find $\delta_i > 0$ small enough and that does not depend on $k$ such that for $k$ large enough:
\begin{equation} \label{minorunif}
\inf_{B_{\check{x}_{i,k}}(\delta_i)} \check{u}_k  \ge \frac{1}{4}| \check{x}_i |^{1 - \frac{n}{2}}. 
\end{equation}
Assume on the contrary that \eqref{uqnonborne} holds. By the definition of $d_k$ in \eqref{defdq}, conditions \eqref{hyposuites} are satisfied with $x_k = x_{i,k}$ and $\rho_k = \frac{1}{8} d_k$ so Lemma \ref{lemma7} in the form of \eqref{asympto} applies and shows that:
\begin{equation} \label{blowupcheckuq}
 \check{u}_k( \check{x}_{i,k}) \to + \infty 
 \end{equation}
and that
\begin{equation} \label{explosionuq}
\check{u}_k( \check{x}_{i,k})^{1 -n + q_k \frac{n-2}{2}} \check{u}_k (x) \to \frac{\lambda_i}{|x- \check{x}_i|^{n-2}} + H_i(x) 
\end{equation}
in $C^1_{loc}\left( B_{\check{x}_i}(\frac{1}{2}) \backslash \{\check{x}_i \} \right)$ as $k \to \infty$, where $\lambda_i >0$ and $H_i$ is a harmonic  function in $B_0(\frac{5}{8})$ satisfying $H_i(\check{x}_i) \le 0$. 
Now we show that these two situations cannot simultaneously happen.  Let $\mathcal{A} \subset \{1 \dots N_k \}$ be a finite collection of subscripts such that there exists a positive $R$ such that for all $i \in \mathcal{A}$ 
\[ d_g(x_{1,k},x_{i,k}) \le R d_k .\]
Assume that there exist $i, j \in \mathcal{A}$, $i < j$, such that $x_{i,k}$ satisfies \eqref{uqborne} and $x_{j,k}$ satisfies \eqref{uqnonborne}. Since $1 -n + q_k \frac{n-2}{2} \to 1$ as $k \to \infty$, \eqref{blowupcheckuq} and the $C^1$-convergence result  \eqref{explosionuq} show that for any positive $r$:
\begin{equation} \label{checkuqconv0}
\check{u}_k \to 0 \textrm{ in } C^1_{loc}\left( B_{ \check{x}_{j}}(2r) \backslash B_{ \check{x}_{j}}(r) \right) .
\end{equation}
By their definition in \eqref{checkxi}, critical points $\check{x}_{l,k}$ satisfy for all $l, m \in \mathcal{A}$, $ l\neq m$, $|\check{x}_{l,k} - \check{x}_{m,k}| \ge 1$. It is then possible to find a connected open set in  $\mathbb{R}^n$, which we will call $U$, $U \subset B_0(R+1)$, that contains $B_{\check{x}_{i,k}}(1/4)$ and $B_{\check{x}_{j,k}}(1/4)$ for any $k$ but that does not contain any other point $\check{x}_{l,k}$ for $l \neq i,j$. Let $0 < r  < \frac{1}{8}$ and consider 
\begin{equation} \label{domaineVrR}
V_{r,R} = U \backslash \left( \overline{B_{\check{x}_{i,k}}(r)} \cup \overline{B_{\check{x}_{j,k}}(r)} \right) \cap \Omega_{r,R+1}.
\end{equation}
By \eqref{petitHarnack} we get
\[ \sup_{V_{r,R}} \check{u}_k \le D_{r,R+1} \inf_{B_{\check{x}_{j,k}}(2r) \backslash \overline{B}_{\check{x}_{j,k}}(r)} \check{u}_k\hskip.1cm , \]
so that using \eqref{checkuqconv0} we obtain that $\check{u}_k$ goes uniformly to $0$ on every annulus centered at $\check{x}_{i,k}$, which contradicts \eqref{minorunif}.
\end{proof}

We are able to conclude the proof of Theorem \ref{Thstabi}. We first rule out the case \eqref{uqborne}. Assume by contradiction that, for any $1 \le i \le N_k$ such that $d_g(x_{1,k},x_{i,k}) = O( d_k )$, \eqref{uqborne} holds.  By \eqref{petitHarnack} and \eqref{uqborne} the sequence $(\check{u}_k)_k$ is uniformly bounded in $L^\infty(B_0(R))$ for all $R > 0$. Using \eqref{eqcheckuq}, \eqref{propcheckfq} and standard elliptic theory there exists $\check{u}$ satisfying
\[ \triangle_{\xi} \check{u} = f(x_1) \check{u}^{2^*-1} \]
and such that
\[ \check{u}_k \to \check{u} \textrm{ in } C^1_{loc}(\RR^n) \]
as $k \to \infty$, where  $x_1 = \lim x_{1,k}$ up to a subsequence. We now show that $f(x_1) > 0$. By \eqref{grossehypo} we therefore assume that $\| u_k \|_{H^1(M)} \le C_0$. In particular the limit function $\check{u}$ belongs to $L^{2^*}(\mathbb{R}^n)$ since for any positive $R$:
\[ \int_{B_0(R)} \check{u}_k^{2^*} d v_{\check{g}_k} =
  O \left( \int_{B_{x_{1,k}}(d_k R)} d_k^{\frac{2q_k}{q_k-2} - n} u_k^{2^*} dv_g \right) = O \left( \int_M u_k^{2^*} dv_g \right) = O(1)\hskip.1cm ,\]
where the last equality holds true because $\frac{2q_k}{q_k-2} - n \ge 0$ and $d_k \to 0$. Passing to the limit \eqref{ptscrit1} with \eqref{propdq} we obtain $\check{u}(0) 
 \ge 1$. But $\check{u}$ cannot be subharmonic, non-zero and belong to $L^{2^*}(\RR^n)$, hence $f(x_1) >0$. Moreover, we know that $\check{u}$ has at least two distinct critical points, $0$ and $\check{x}_2$, since by assumption \eqref{propdq} there holds $|\check{x}_{2,k}| = 1$ for all $k$. This contradicts the classification result in Caffarelli-Gidas-Spruck \cite{CaGiSp} and thus contradicts \eqref{uqborne}. Hence for any $1 \le i \le N_k$ such that $d_g(x_{1,k},x_{i,k}) = O (d_k)$ we have
\[ \check{u}_k(\check{x}_{i,k}) \to + \infty \]
as $k  \to \infty$. Proceeding as in the proof of Lemma \ref{lemma10}, using the Harnack inequality \eqref{petitHarnack} it is easily seen that $\check{u}_k$ blows up at the same speed at every concentration point at a finite distance from $0$: for any $1 \le i \le N_{k}$ such that $d_g(x_{1,k},x_{i,k}) = O( d_k)$ there exist $\mu_i>0$ such that up to a subsequence
\begin{equation} \label{explosionrelative}
\frac{\check{u}_k(\check{x}_{i,k})}{\check{u}_k(0)} \to \mu_i 
\end{equation}
as $k \to \infty$. Using \eqref{explosionrelative}, \eqref{explosionuq} shows that if $d_g(x_{1,k}, x_{i,k}) = O(d_k)$, then
\begin{equation} \label{convautourxi}
 \check{u}_k(0)^{1 -n + q_k \frac{n-2}{2}} \check{u}_k (x) \to \frac{\lambda_i}{\mu_i |x- \check{x}_i|^{n-2}} + \mu_i^{-1} H_i(x)  
 \end{equation}
in $C^1_{loc}(B_{\check{x}_i} (\frac{1}{2}) \backslash \{ \check{x}_i \})$ as $k \to \infty$. Let $R>0$. We know that $(N_{R,k})_k$ and $(N_{2R,k})_k$ are bounded and we can thus assume they are constants $N_R$, $N_{2R}$. Using Harnack's inequality \eqref{petitHarnack} and \eqref{convautourxi} it is easily seen that 
\begin{equation} \label{allurecheckG}
 \check{u}_k(0)^{1-n+q_k \frac{n-2}{2}} \check{u}_k \to \check{G} \textrm{ in } C^1_{loc}(B_0(R) \backslash \{ \check{x}_i \}_{i=1, \cdots N_{2R}}) \hskip.1cm ,
 \end{equation}
where, by definition of $N_R$ as in \eqref{defNrk}, we can write that
\begin{equation} \label{checkG}
\check{G}(x) = \sum_{i=1}^{N_R} \frac{\lambda_i}{\mu_i |x - \check{x}_i|^{n-2}} + \check{H}(x)\hskip.1cm ,
\end{equation}
and where $\check{H}$ is a harmonic function in $\overline{B_0(R)}$. We write 
\[ \check{G}(x) = \frac{\lambda_1}{|x|^{n-2}}  + \left( \sum_{i=2}^{N_R} \frac{\lambda_i}{\mu_i |x - \check{x}_i|^{n-2}} + \check{H}(x) \right)  \]
and since $\check{u}_k(0) \to \infty$ as $k\to\infty$ we apply Lemma \ref{lemma7} with $x_k = x_{1,k}$ and $\rho_k = \frac{1}{16} d_k$ to obtain
\begin{equation} \label{condH1}
 \sum_{i=2}^{N_R} \frac{\lambda_i}{\mu_i |\check{x}_i|^{n-2}} + \check{H}(0) \le 0. 
\end{equation} 
Independently, let $z_0$ be some point satisfying $|z_0| = R$ and $|z_0-\check{x}_2| \ge R-1$ and different from any $\check{x}_i$, $3 \le i \le N_{2R}$. Let $U_0$ be some connected open set in $\overline{B_0(R)}$ that contains $0$ in its interior and $z_0$ on its boundary and avoids any other $\check{x}_i$. Then $\check{G} - \lambda_1 |x|^{2-n} - \frac{\lambda_2}{\mu_2} |x-\check{x}_2|^{2-n}$ is harmonic in $U_0$ , so by the maximum principle: 
\begin{equation} \label{ppemaxcheckG}
\begin{aligned} 
\left( \check{G} - \frac{\lambda_1}{|x|^{n-2}} - \frac{\lambda_2}{\mu_2 |x-\check{x}_2|^{n-2}} \right)(0) & \ge \left( \check{G} - \frac{\lambda_1}{|x|^{n-2}} - \frac{\lambda_2}{\mu_2 |x- \check{x}_2|^{n-2}} \right)(z_0) \\
& \ge - \frac{\lambda_1}{R^{n-2}} - \frac{\lambda_2}{\mu_2 (R-1)^{n-2}} 
\end{aligned}
\end{equation}
since $\check{G}(x) \ge0$ in $B_0(R) \backslash \{ 0\}$. With \eqref{checkG} and \eqref{ppemaxcheckG} we obtain in the end, since $|\check{x}_2| = 1$: 
\[  \sum_{i = 2}^{N_R} \frac{\lambda_i}{\mu_i |\check{x}_i|^{n-2}} + \check{H}(0) \ge \frac{\lambda_2}{\mu_2} - \frac{\lambda_1}{R^{n-2}} - \frac{\lambda_2}{\mu_2 (R-1)^{n-2}}
\hskip.1cm ,  \]
which contradicts \eqref{condH1} up to choosing $R$ large enough. This shows that \eqref{blowup} can never happen and concludes the proof of Theorem \ref{Thstabi}.

\section{A general mountain-pass situation.} \label{mountainpass}

In this section we prove, as an application of Theorem \ref{Thstabi}, a result that states that equation \eqref{EL} has at least two solutions as soon as it has a mountain-pass structure. Starting with \eqref{EL}, we introduce the following equations for all $2 \le q \le  2^*$:
\begin{equation} \label{eq:einlichq} \tag{$EL_q$}
\triangle_g u + h u  = f u^{q-1} + \frac{a}{u^{q+1}}
\end{equation}
and the associated energy functional, defined on $H^1(M)$ whenever the following expression makes sense (for instance for positive functions):
\begin{equation} \label{fonctionnelleq}
I^q(u) = \frac{1}{2} \int_M \left( |\nabla u|_g^2 + hu^2 \right) dv_g - \frac{1}{q} \int_M f (u^+)^{q} dv_g + \frac{1}{q} \int_M \frac{a}{(u^+)^{q}}dv_g. 
\end{equation}
Equation \eqref{eq:einlichq} is called subcritical if $q < 2^*$ and critical if $q = 2^*$. For any $2 \le q \le 2^*$ we call $S_{h,q} = S_{h,q}(M,g)$ the smallest positive constant in the embedding $H^1(M) \subset L^{q}(M)$ for the $H^1_h$-norm, i.e. satisfying:
\begin{equation} \label{defSh}
\Vert u \Vert_{L^{q}} \le S_{h,q}^{\frac{1}{q}} \Vert u \Vert_{H^1_h}\hskip.1cm ,
\end{equation}
where the $H^1_h$-norm is as in \eqref{Nh}. We have $S_{h,2^*} = S_h$, where $S_h$ is as in \eqref{Sh}. Following Hebey-Pacard-Pollack \cite{HePaPo}, to get rid of the negative exponent we consider  for any $\ve > 0$ the perturbed functional on $H^1(M)$:
\begin{equation} \label{Iqeps}
I_\ve^q  = \frac{1}{2} \int_M \left( |\nabla u|_g^2 + hu^2 \right) dv_g - \frac{1}{q} \int_M f (u^+)^{q} dv_g + \frac{1}{q} \int_M \frac{a}{\left( \ve + (u^+)^2 \right)^{q/2}}dv_g\hskip.1cm ,
\end{equation}
where $2 \le q \le 2^*$. Our result states that, in low dimensions, each time the critical equation \eqref{EL} has a mountain-pass structure each equation \eqref{eq:einlichq} with $2 \le q \le 2^*$ actually has two distinct positive solutions. We state it as follows.

\begin{theo} \label{Th01}
Let $(M,g)$ be a $n$-dimensional closed Riemannian manifold with $3 \le n \le 5$ and $h, f$ and $a$ be smooth functions in $M$ such that $\triangle_g+h$ is coercive, $a \ge 0, a \not \equiv 0$ and $\max_M f >0$. We assume that there exist two smooth positive functions $u_0$ and $u_1$ and two real numbers $\eta >0$ and $\rho >0$ such that $\Vert u_1 - u_0\Vh < \rho$, 
\begin{equation} \label{inf2}
I^{2^*}(u_1)< \eta\hskip.1cm ,
\end{equation}
where $I^{2^*}$ is as in \eqref{fonctionnelleq}, and
\begin{equation} \label{inf}
\inf_S I^q_\ve(u) > \eta
\end{equation} 
for all $\ve$ small enough and for all $ q $ close enough to $2^*$, where we have let $S$ be the boundary set $S = \{ u \in H^1(M),  \Vert u - u_0\Vert_{H^1_h} = \rho \} $. Then for any $q$ close enough to $2^*$ equation \eqref{eq:einlichq} admits two distinct smooth positive solutions.
\end{theo}

\begin{proof}
We first consider the case $q < 2^*$ and then prove the $q = 2^*$ case with a stability argument. We thus assume that $2 \le q < 2^*$. The first solution is obtained as a limit of local minima of $I_\ve^q$. The second one is obtained as a mountain-pass solution. On the closed ball 
\begin{equation} \label{balle}
B = \{ u \in H^1(M), \Vert u - u_0\Vh \le \rho \}
\end{equation}
the functional $I_\ve^q$ is uniformly (in $\ve$) bounded from below. Indeed, since $u_1$ is positive by Lebesgue's dominated convergence theorem we have: 
\[ \lim_{q \to 2^*} \lim_{\ve \to 0} I_\ve^q(u_1)  =  I^{2^*}(u_1)\hskip.1cm ,\] 
so that using \eqref{inf2} there holds then for $\ve$ small enough and for $q$ close to $2^*$,
\begin{equation} \label{borneunif1}
 -\frac{1}{q}\max_M |f| S_{h,q} (\rho + \Vert u_0 \Vert_{H^1_h})^q \le \inf_B I_q^\ve \le \eta\hskip.1cm ,
\end{equation}
where $S_{h,q}$ is as in \eqref{defSh} and $I_q^\ve$ as in \eqref{Iqeps}. Using Ekeland's variational principle we can let $(u_k)_k$ be a Palais-Smale minimising sequence for $I_\ve^q$ such that $I_\ve^q(u_k) \to \inf_B I_\ve^q$ as $k\to\infty$ (see for instance Struwe \cite{Stru}, Chapter $1$, Corollary $5.3$ for the proof of this statement). The sequence $(u_k)_k$ thus satisfies
\begin{equation} \label{approx1}
\frac{1}{2} \int_M \left( |\nabla u_k|_g^2 + h u_k^2 \right) dv(g) - \frac{1}{q} \int_M f (u_k^+)^{q} dv(g) + \frac{1}{q} \int_M \frac{a}{(\varepsilon + (u_k^+)^2)^{\frac{q}{2}}} = \inf_B I_\ve^q + o(1) 
\end{equation}
and, for all $\Phi_k \in H^1(M)$,
\begin{equation} \label{approx2}
\begin{aligned}
\int_M \langle \nabla u_k, \nabla \Phi_k \rangle_g + h u_k \Phi_k dv(g) - & \int_M f (u_k^+)^{q-1} \Phi_k dv(g) \\
 & - \int_M \frac{a u_k^+}{(\varepsilon + (u_k^+)^2)^{\frac{q}{2}+1}} \Phi_k dv(g) = o \left( \Vert \Phi_k \Vert_{H^1} \right) . \\
\end{aligned}
\end{equation}
Choosing $\Phi_k = u_k$ in \eqref{approx2} and combining with \eqref{approx1} yields, since $a \ge 0$:
\begin{equation} \label{approx3}
\left( \frac{1}{2} - \frac{1}{q} \right) \int_M f (u_k^+)^{q} dv(g) \leqslant  \inf_B I_\ve^q + o \left( ||u_k||_{H^1_h} \right) + o(1).
\end{equation} 
Using \eqref{approx3} in \eqref{approx1} we get for $k$ large enough:
\begin{equation} \label{approximp}
\begin{aligned}
\Vert u_k \Vh^2 \le \frac{4q}{q-2} \inf_B I_\ve^q + o(1) 
\\
\end{aligned}
\end{equation}
and in particular $ \inf_{B} I_\ve^q \ge 0$ .
There exists then $u_{\ve,q} \in H^1(M)$ such that, up to a subsequence, $u_k \rightharpoonup u_{\ve,q}$ in $H^1(M)$.  By standard integration theory $u_{\ve,q}$ satisfies weakly
\begin{equation} \label{Th13}
 \triangle_g u_{\varepsilon,q} + h u_{\varepsilon,q} = f (u_{\varepsilon,q}^+)^{q-1} + \frac{a (u_{\varepsilon,q}^+)}{(\varepsilon + (u_{\varepsilon,q}^+)^2)^{\frac{q}{2}+1}} .
\end{equation}
Multiplying \eqref{Th13} by $u_{\ve,q}^{-}$ and integrating it is easily seen that $u_{\ve,q}$ is nonnegative almost everywhere. Since $a (\varepsilon + (u_{\varepsilon,q}^+)^2)^{- \frac{q}{2}-1}$ belongs to $L^\infty(M)$ standard bootstrap arguments show that $u_{\ve,q}$ is smooth in $M$. 
Now we obtain a uniform bound from below on $u_{\ve,q}$ and show that $\inf_M u_{\ve,q}$ does not converge to $0$ as $\ve \to 0$. We consider for any positive $\delta$ the unique functions $\psi_\delta$ and $\psi_0$ solving 
\begin{equation} \label{astucedelta}
\begin{aligned}
&  \triangle_g \psi_\delta + h \psi_\delta = a - \delta f^-  -\delta \\
 & \triangle_g \psi_0 + h \psi_0 = a  .
  \end{aligned}
\end{equation}
Since $a \not \equiv 0$, $\psi_0 > 0$ in $M$. By standard elliptic theory $\psi_{\delta} \to \psi_0 $ in $C^0(M)$ when $\delta \to 0$ and thus $\psi_{\delta_0}$ is positive for some $\delta_0 > 0$ small enough. Let $\ve_0$ small enough so that $\inf_M \psi_{\delta_0} \ve_0^{- 2^*-2} > 1$ and let $t_0>0$ small enough such that for any $q$ close enough to $2^*$ 
\[ \frac{ \psi_{\delta_0}(x)}{\left( \ve_0^2  +  t_0^2 \psi_{\delta_0}(x)^2 \right)^{\frac{q}{2}+1}} > 1 \]
for any $x \in M$. Then for any $\ve \le \ve_0$ and $t \le t_0$, $\theta_t = t \psi_{\delta_0}$ satisfies
\begin{equation} \label{soussoldelta}
\triangle_g \theta_t + h \theta_t < f \theta_t^{q-1} + \frac{a \theta_t}{\left( \ve + \theta_t^2 \right)^{\frac{q}{2}+1}} .
\end{equation}
Now we claim that there exists some $t > 0$ small enough such that for any $\ve$ small enough and any $q$ close enough to $2^*$ and for any smooth positive $\vp $ solution of \eqref{Th13}  there holds 
\begin{equation} \label{minunifveq} 
\vp > \theta_t
\end{equation} 
in $M$. This shows in particular that $u_{\ve,q}$ is uniformly bounded from below in $\ve$ and $q$. We prove the claim by contradiction and assume that for any positive $t$ there exists $\ve_t > 0$, $2 \le q_t < 2^*$, $\vp$ a solution of \eqref{Th13} with $\ve = \ve_t$ and $q = q_t$ and $x_t \in M$ such that $\vp(x_t)  \le \theta_t(x_t)$. Then, for some $\tilde t \in (0,t)$, and some $\tilde x_t \in M$, 
\begin{equation} \label{infuveq}
 1 = \inf_M \frac{\vp}{\theta_{\tilde t}} = \frac{ \vp (\tilde x_t)}{\theta_{\tilde t}(\tilde x_t)}~.
 \end{equation}
In particular, with \eqref{infuveq} we obtain that
\[\theta_{\tilde t}(\tilde x_t) = \vp (\tilde x_t)~~\hbox{and}~\triangle_g \vp (\tilde x_t) \le \triangle_g \theta_{\tilde t} (\tilde x_t) \]
which is impossible since $\theta_{\tilde t}$ is a strict subsolution of \eqref{Th13} by \eqref{soussoldelta}. This proves \eqref{minunifveq} which shows that $u_{\ve,q}$ is positive for any $\ve$ small enough and any $q$ close to $2^*$. Finally since $q < 2^*$, $u_k \to u_{\ve,q}$ in $L^q(M)$. Since $u_{\ve,q}$ solves  \eqref{Th13}, a straightforward computation using \eqref{approx2} with $ \Phi = u_k$ yields:
\begin{equation} \label{grad}
\int_M |\nabla u_k - \nabla u_{\varepsilon,q}|_g^2 dv_g   = o(1) 
\end{equation}
as $k \to \infty$ and so $u_k \to u_{\ve,q}$ in $H^1(M)$.  As a consequence $I_\ve^q(u_{\ve,q}) = \inf_B I_\ve^q$. In particular by definition of $I_\ve^q$ in \eqref{Iqeps} there holds $ I_\ve^q(u_1) \le I^q(u_1)$, where $I^q$ is as in \eqref{fonctionnelleq}. Since $u_1 \in B$, with $B$ as in \eqref{balle}, using  \eqref{inf2} there holds, for $q$ close enough to $2^*$:
\begin{equation} \label{Th14}
I_\ve^q(u_{\ve,q}) \le I^q(u_1)< \eta.
\end{equation}
Now we let $\ve \to 0$. By \eqref{approximp} and  \eqref{borneunif1} there exists a nonnegative $u_q \in H^1(M)$ such that $u_{\ve,q}$ converges, up to a subsequence, weakly to $u_q$  in $H^1(M)$. With \eqref{minunifveq} and standard integration theory we can pass to the limit in \eqref{Th13} and get that $u_q$ is a smooth solution of \eqref{eq:einlichq}, once again positive by \eqref{minunifveq}. 
Here again, since $q < 2^*$, a similar computation to the one in \eqref{grad} shows that 
\[ \int_M |\nabla u_{\ve,q} - \nabla u_q|^2 dv_g = o(1) \]
as $\ve \to 0$, so $u_{\ve,q} \to u_q$ strongly in $H^1(M)$. In particular \eqref{minunifveq} and Lebesgue's dominated convergence theorem show that $I^q_\ve(u_{\ve,q}) \to I^q(u_q)$, where $I^q$ is as in \eqref{fonctionnelleq}. By \eqref{Th14} this gives, for $q$ close enough to $2^*$:
\begin{equation} \label{Th15}
I^q(u_{q})  \le I^q(u_1) < \eta.
\end{equation}
Now we take into account the mountain-pass structure of $I^q_\ve$ in order to construct a second solution. Since $\max_M f >0$ it is possible to find a smooth positive function $\psi $ in $M$ such that 
\begin{equation} \label{intf}
\int_M f \psi^{2^*} dv_g > 0 .
\end{equation}
Then by \eqref{fonctionnelleq} and \eqref{intf}, $I^{2^*}(t \psi) \to - \infty$ as $t$ goes to infinity. We can thus pick a large real number $T$ such that, for any $q$ close enough to $2^*$ and for any positive $\ve$ small enough,  
\begin{equation} \label{hypocolint}
I^q_\ve(T \psi) < \eta
\end{equation} and also $\Vert T \psi - u_0 \Vert_{H^1_h} > \rho$, where $\eta, u_0$ and $\rho$ are as in the statement of Theorem \ref{Th01}. On the other side by assumption \eqref{inf} we get for any $q$ close enough to $2^*$ and for any positive $\ve$ small enough
\begin{equation} \label{hypocolinf}
I^q_\ve(u) > \eta 
\end{equation}
for all $u \in S= \{ u \in H^1(M), \Vert u - u_0\Vh = \rho \}$. Thus with \eqref{Th14}, \eqref{hypocolint} and \eqref{hypocolinf}, for any $q$ close enough to $2^*$, and for any $\ve$ small enough, we can apply the mountain-pass lemma as stated in Rabinowitz \cite{Rab} that provides us with a Palais-Smale sequence $(v_k)_k$ such that
\begin{equation} \label{Th16}
I^q_\ve(v_k) \to c_{\ve,q} 
\end{equation}
as $k \to +\infty$, where we have set
\begin{equation} \label{defcq}
c_{\ve,q} = \inf_{h \in \Gamma} \max_{u\in h([0;1])} I^q_\ve(u) 
\end{equation}
and where $\Gamma$ is the set of paths $h: [0,1] \to H^1(M)$ such that $h(0) = u_{\ve,q}$ and  $h(1) = T \psi$. Using \eqref{inf} and considering as a special $h$ a parametrization of the segment $[u_{\ve,q}; T \psi]$ it is easily seen, using \eqref{approximp} and \eqref{minunifveq} that there exists a positive $C_0$ such that 
\begin{equation} \label{majorcveq}
\eta \le c_{\ve,q} \le C_0 
\end{equation} 
for $\ve$ small enough and for $q$ close to $2^*$. Therefore the same computation that led to \eqref{approximp} 
works again and shows that
\begin{equation} \label{approximpv}
\| v_k \|_{H^1_h} \le \frac{4q}{q-2} C_0 + 1
\end{equation}
with $C_0$ as in \eqref{majorcveq}. With \eqref{approximpv} and standard integration theory, mimicking the proof of \eqref{grad} one shows that the sequence $(v_k)_k$ converges strongly in $H^1(M)$ to some nonnegative function $v_{\ve,q}$ that solves \eqref{Th13}. Passing in \eqref{hypocolinf} to the limit as $k \to \infty$ shows that
\begin{equation} \label{infvveq}
I^q_\ve(v_{\ve,q}) = c_{\ve,q} \ge \eta
\end{equation}
by \eqref{majorcveq}. Also with \eqref{minunifveq} we get that $v_{\ve,q}$ is uniformly bounded from below in $\ve$ and in $q$. Thus repeating the whole argument once again it is easily seen that $v_{\ve,q}$ converges strongly in $H^1(M)$ as $\ve \to 0$ to some smooth positive solution $v_q$ of \eqref{eq:einlichq}. 
With \eqref{minunifveq} we can pass in \eqref{infvveq} to the limit as $\ve \to 0$ and obtain: 
\begin{equation} \label{valeur}
I^q(v_q) \ge \eta.
\end{equation}
In particular, \eqref{Th15} and \eqref{valeur} show that $u_q \not = v_q$. This shows Theorem \ref{Th01} for all $2 \le q < 2^*$ and in any dimension. We now conclude the proof and show that there still exist two different solutions of the critical equation in low dimensions, that is when $3 \le n \le 5$. By \eqref{approximp} and \eqref{approximpv} the sequences $(u_q)_q$ and $(v_q)_q$ are bounded in $H^1(M)$. Since $\triangle_g + h$ is coercive and $\max_M f >0$ we can let $(q_k)_k$ be some sequence converging to $2^*$ and apply Theorem \ref{Thstabi} with $a_k \equiv a $ to the sequences $(u_{q_k})_k$ and $(v_{q_k})_k$. There thus exist $u$ and $v$ two smooth positive functions that solve \eqref{EL} such that $u_{q_k} \to u$ and $v_{q_k} \to v$ in $C^{1, \alpha}(M)$ for all $0 < \alpha < 1$. Passing \eqref{Th15} and \eqref{valeur} to the limit as $k \to \infty$ finally gives:
\[ I^{2^*}(u) \le I^{2^*}(u_1) < \eta \le I^{2^*}(v) \]
which shows that $u$ and $v$ are distinct.
\end{proof}
Note that the construction of $u$ and $v$ is consistent with the computation of the degree of \eqref{EL} we performed in \eqref{calculdegre}. Indeed, $u$ is obtained as a limit of local minima, hence of solutions of index $0$ while $v$ is a mountain-pass solution, hence (generically) of index $+1$.

 \section{A minimal solution of \eqref{EL}.} \label{solumini}

We now investigate more precisely the influence of the parameter $a$ and consider the following equation: 
\begin{equation} \label{ELa} \tag{$EL_a$}
\triangle_g u + h u  = f u^{2^*-1} + \frac{a}{u^{2^*+1}}\hskip.1cm ,
\end{equation}
where $h, f$ and $a$ are smooth functions in $M$, $\triangle_g + h$ is coercive, $\max_M f >0$ and $a$ is nonnegative and nonzero. Using the sub and super solution method we show that each time equation \eqref{ELa} has a smooth positive solution it has a smallest solution for the $L^\infty$-norm:

\begin{prop} \label{solmin}
Let $(M,g)$ be a closed Riemannian manifold of dimension $n \ge 3$ and $a$ be a non zero smooth function in $M$. Let $h$ and $f$ be smooth functions in $M$ such that $\triangle_g + h$ is coercive and $\max_M f > 0$ and $a  \ge 0$. Assume that  \eqref{ELa} has a smooth positive solution. Then there exists a smooth positive function $\vp(a)$ solving \eqref{ELa} such that for any other solution $\vp$ of \eqref{ELa} with $\vp \not \equiv \vp(a)$ there holds $\vp > \vp(a)$. Moreover $\vp(a)$ is stable, in the sense that for any $\psi \in H^1(M)$ there holds
\[ \int_M |\nabla \psi|^2 + \left [h - (2^*-1)f \vp(a)^{2^*-2} + (2^*+1) \frac{a}{ \vp(a)^{2^*+2}} \right] \psi^2 dv_g \ge 0 ,\]
and the mapping $a \mapsto \vp(a)$ in nondecreasing in the following sense: if $a_1 \le a_2$ in $M$, provided $\vp(a_1)$ and $\vp(a_2)$ exist, there holds $\vp(a_1) \le \vp(a_2)$.
\end{prop}

\begin{proof}
Let $a \ge 0$ be a nonzero smooth function such that \eqref{ELa} has a solution. Mimicking the proof of \eqref{minunifveq} we start proving that there exists a positive number that bounds from below all the solutions of \eqref{ELa}. As in \eqref{minunifveq}, notice that there always exist sub-solutions of \eqref{ELa} as small as we want. Indeed, for any $\delta \ge 0$ we let $u_\delta$ be the unique solution of
\begin{equation} \label{udelta}
 \triangle_g u_\delta + h u_\delta = a - \delta f^-  -\delta\hskip.1cm ,
 \end{equation}
where $f^- = - \min(f,0)$. Since $a$ is nonnegative and nonzero, the maximum principle shows that $u_0 > 0$ in $M$. By standard elliptic theory $\Vert u_\delta - u_0 \Vert_\infty \to 0$ as $\delta$ goes to $0$ so for some $\delta_0 > 0 $ small enough we have $u_{\delta_0} >0$. Then for $\ve$ small enough, 
\begin{equation} \label{soussol}
v_\ve = \ve u_{\delta_0}
\end{equation}
is a strict sub-solution of \eqref{ELa} since, by \eqref{udelta},
\[ \triangle_g v_\ve + h v_\ve = \ve a - \ve \delta_0 f^-  - \ve \delta_0< \frac{a}{v_\ve^{2^*+1}} + f v_\ve^{2^*-1}. \]
Now we claim that there exists some $\ve_0 > 0$ such that for any positive solution $\vp$ of \eqref{ELa} there holds 
\begin{equation} \label{minunif}
\vp > v_{\ve_0}
\end{equation} 
in $M$, where $v_{\ve_0}$ is as in \eqref{soussol}. We prove the claim by contradiction and assume that 
there exists $\vp_\ve$ solution of \eqref{ELa}, and $x_\ve \in M$, such that $\vp_\ve(x_\ve) \le v_\ve(x_\ve)$ for all 
$\ve > 0$. Then, for some $\tilde\ve \in (0,\ve)$, and some $\tilde x_\ve \in M$, 
\[ 1 = \inf_M \frac{\vp_\ve}{v_{\tilde\ve}} = \frac{\vp_\ve(\tilde x_\ve)}{v_{\tilde\ve}(\tilde x_\ve)}~.\]
In particular, we obtain that
\[v_{\tilde\ve}(\tilde x_\ve) = \vp_\ve(\tilde x_\ve)~~\hbox{and}~\triangle_g \vp_\ve (\tilde x_\ve) \le \triangle_g v_{\tilde\ve} (\tilde x_\ve) \]
which is impossible since $v_{\tilde\ve}$ is a strict subsolution of \eqref{ELa}. 
Now we prove the existence of a minimal solution of \eqref{ELa}. We follow here the arguments in Sattinger \cite{Sattinger}. For $x \in M$ and $u >0$ we let 
\[F(x,u) = f(x) u(x)^{2^*-1} + \frac{a(x)}{u(x)^{2^*+1}} - h(x) u(x).\] 
Let $\psi$ be a solution of \eqref{ELa} and let $w$ be a strict subsolution of \eqref{ELa} which is less than any positive solution of \eqref{ELa}. We proved the existence of such a $w$ 
in \eqref{minunif}. Also we let $K > 0$ be large enough such that for any $x\in M$, and any $\min_M w  \le u \le \max_M \psi$,
\begin{equation} \label{propK}
\begin{aligned}
 F(x,u) + Ku \ge 0 \textrm{ and } \frac{\partial F}{\partial u}(x,u) + K \ge 0.
 \end{aligned}
 \end{equation} 
 For any smooth positive function $u$, we define $Tu$ as the unique solution of
\begin{equation} \label{defiT}
\triangle_g Tu + K Tu = F(\cdot,u) + Ku .
\end{equation}
As a first remark, for any two positive functions $u$ and $v$ in the range 
\[ \min_M w \le u,v \le \max_M \psi \]
we have: 
\[ \big( \triangle_g + K \big)(Tu - Tv)(x) = F(x,u) - F(x,v) + K \big( u(x) - v(x) \big).\]
Then, by the strong maximum principle, we obtain that
\begin{equation} \label{croissance}
Tu < Tv \textrm{ as long as } u \le v ~\textrm{and } u \not \equiv v.
\end{equation} 
The iterative sub and super solution method applied in the range $w \le \vp \le \psi$ and starting from the strict sub-solution $w$ provides a sequence $v_n = T^n w$ which is non decreasing by the maximum principle and converges to a fixed-point of $T$, that is to say a solution of \eqref{ELa} (see \cite{Sattinger} for more details). We shall call this solution $\vp(a)$:
\begin{equation} \label{defsol}
\vp(a) = \lim_{n \to \infty} T^n w.
\end{equation}
By standard elliptic theory, $\vp(a)$ is smooth. Note in passing that all the above arguments still work if we only assume that $a$ is continuous, but in this case $\vp(a)$ constructed as in \eqref{defsol} will only be of class $C^{1,\alpha}$ for any $0 < \alpha < 1$. 
Now we show that $\vp(a)$ does not depend on $\psi$ and on $w$. First, $\vp(a)$ as in \eqref{defsol} does not depend on $\psi$. We let $\psi_1$ and $\psi_2$ be two solutions of \eqref{ELa}. We let $K_i$, $i=1,2$ be positive constants satisfying \eqref{propK} in $[ \min_M w ; \max_M \psi_i]$, $T_i$ be the operator defined as in \eqref{defiT} and $\vp_i$ the associated solution as in \eqref{defsol}. Since $(T_1^n w)$ is non decreasing there holds $\vp_1 \ge w$. If we assume for instance that $\max_M \psi_1 \le \max_M \psi_2$ then $\vp_1 \in [ \min_M w ; \max_M \psi_2]$ and thus, by \eqref{defsol} and the maximum principle there holds $\vp_2 \le \vp_1$ 
since $T_2(\vp_1) = \vp_1$. But then $\vp_2$ is a solution of \eqref{ELa} with $\min_M w \le \vp_2 \le \max_M \psi_1$ and thus, once again by the maximum principle, $\vp_1 \le \vp_2$. This proves that $\vp(a)$ does not depend on $\psi$. 
Now we prove that $\vp(a)$ does not depend on the strict subsolution $w$, provided that $w$ is less than any positive solution of \eqref{ELa}. Indeed, for any $\psi$ solution of \eqref{ELa}, if $w_1$ and $w_2$ are 
two such subsolutions, and $\vp_1$ and $\vp_2$ are the associated solutions as in \eqref{defsol}, there holds $w_1 \le \vp_2$ and $w_2 \le \vp_1$. We conclude once again with the maximum principle that 
$\vp_1 \le \vp_2$ and $\vp_2 \le \vp_1$. By the definition of $\vp(a)$ in \eqref{defsol}, and what we just proved, 
for any $\psi$ solution of \eqref{ELa} there holds that $w < \vp(a) \le \psi$, where $w$ is a subsolution that is less than any solution of \eqref{ELa}. With \eqref{croissance} we obtain the desired property:
\begin{equation} \label{minimal}
\vp(a) < \psi \textrm{ or } \vp(a) \equiv \psi.
 \end{equation}
The stability of $\vp(a)$ is a consequence of the minimality of $\vp(a)$. We denote by $\lambda_0$ the first eigenvalue of the linearized operator of equation \eqref{ELa} at $\vp(a)$. The stability of $\vp(a)$ as stated in Proposition \ref{solmin} amounts to say that $\lambda_0 \ge 0$. Assume by contradiction that $\lambda_0 < 0$ and denote by  $\psi_0$ the associated positive eigenvector. Let $w$ be a subsolution that is less than any solution of \eqref{ELa}. Let $\vp_\delta = \vp(a) - \delta \psi_0$ for any positive $\delta$. For $\delta > 0$ small enough one has
\[ w <  \vp_\delta < \vp(a) \]
and a straightforward calculation shows that
\[ \triangle_g  \vp_\delta +h \vp_\delta - f \vp_\delta^{2^*-1} - \frac{a}{\vp_\delta^{2^*+1}} =  - \delta \lambda_0 \psi_0 + o(\delta) > 0\] 
so that $\vp_\delta$ is a strict supersolution of \eqref{ELa} satisfying $w < \vp_\delta < \vp(a)$ for $\delta$ small enough. By the iterative sub and super solution method we then get a solution
$\psi$ of $(EL_a)$ such that $w < \psi < \varphi_\delta$, and this is in contradiction with \eqref{minimal}. 
Finally, if $a_1 \le a_2$ are nonnegative nonzero functions on $M$, $\vp(a_2)$ is a super solution of equation \eqref{ELa} with $a = a_1$. By the minimality of $\vp(a_1)$ we then have $\vp(a_1) \le \vp(a_2)$.
 \end{proof}

\section{Multiplicity of solutions of \eqref{EL}} \label{deuxsol}

Let $(M,g)$ be a closed $n$-dimensional Riemannian manifold, $3 \le n \le 5$, and $h,f,a$ be smooth functions in $M$ with $a$ nonzero and nonnegative, $\triangle_g + h$ coercive and $\max_M f > 0$. Using Theorem \ref{Th01} and Proposition \ref{solmin} we conclude in this section the proof of Theorem \ref{Th1}. Remember that we are investigating the number of solutions of the following equation:
\begin{equation} \label{eq:einlicht} \tag{$EL_{\theta}$}
\triangle_g u + hu = f u^{2^*-1} + \frac{\theta a }{u^{2^*+1}},
\end{equation}
according to the value of the positive parameter $\theta$. 

\subsection{Two solutions when $a$ is small.}
We now prove the existence of $\theta_1$ as in the statement of Theorem \ref{Th1}. Following Hebey-Pacard-Pollack \cite{HePaPo} we prove that if there exists a positive function $\varphi$ in $M$ satisfying $\| \vp \|_{H^1_h} = 1$ and
\begin{equation} \label{eq:preuveexis}
\int_M \frac{a}{\varphi^{2^*}} dv_g \leqslant \frac{C(n)}{(S_h \max_M |f|)^{n-1}}\hskip.1cm ,
\end{equation}
where 
\begin{equation} \label{C(n)}
C(n) = \frac{1}{n-2} \frac{1}{\left( 2(n-1)\right)^{\frac{2^*}{2}}}\hskip.1cm ,
\end{equation}
then \eqref{EL} has at least two smooth positive solutions. Note that for the sake of simplicity in \eqref{eq:preuveexis} we have let $S_h = S_{h, 2^*}$, where $S_{h,2^*}$ is as in \eqref{defSh}. To do this, we prove that if \eqref{eq:preuveexis} is satisfied so are the assumptions of Theorem \ref{Th01}. We introduce for any $2 \le q \le 2^*$ the function
\[ \Phi_q(t)  = \frac{1}{2}t^2 - \frac{\max_M |f|}{q}S_{h,q} t^{q} \]
defined on $\mathbb{R}^+$ which attains its maximum at 
\begin{equation} \label{t0}
t_{0,q}   = \left(\frac{1}{S_{h,q} \max_M |f|} \right)^{\frac{1}{q-2}}\hskip.1cm ,
\end{equation}
where $S_{h,q}$ is as in \eqref{defSh}. The value of $\Phi_q$ at its maximum is, for $2 \le q \le 2^*$, given by
\begin{equation} \label{ineqq}
\Phi_q(t_{0,q}) =  \left( \frac{1}{2} - \frac{1}{q} \right) \frac{1}{\left(S_{h,q} \max_M |f|\right)^{\frac{2}{q-2}}} .
\end{equation}
In particular, since $S_{h,q} \to S_{h,2^*}$ as $q \to 2^*$:
\begin{equation} \label{limphiq}
\Phi_q(t_{0,q}) \to \Phi_{2^*}(t_{0,2^*}) > 0
\end{equation}
as $q \to 2^*$. We shall write for the sake of simplicity $\Phi(t_0) = \Phi_{2^*}(t_{0,2^*})$ in the following, where we have thanks to \eqref{ineqq}:
\begin{equation} \label{phi0}
\Phi(t_0) = \frac{1}{n} \left( S_{h} \max_M |f| \right)^{- \frac{n-2}{2}}.
\end{equation}
By their definition in \eqref{Iqeps}, the perturbed functionals $I_\ve^q$ satisfy for any $u \in H^1(M)$:
\begin{equation} \label{ineq1}
\begin{aligned}
I_\ve^q(u) \ge  \Phi_q(||u||_{H_1^h}) .
 \end{aligned}
\end{equation}
In particular, from \eqref{ineq1} there holds for any $2 \le q \le 2^*$, any $\ve$ and any $\Vert u \Vh = t_{0,q}$:
\begin{equation} \label{ineq2}
I^q_\ve(u) \ge \Phi_q(t_{0,q}) .
\end{equation}
We now let
\begin{equation} \label{t1}
t_1 = \left( 2(n-1)\right)^{- \frac{1}{2}} t_0.
\end{equation}
As one can easily check thanks to \eqref{t0}, \eqref{phi0} and \eqref{t1} there holds:
\begin{equation} \label{valt1}
 \frac{1}{2}t_1^2 + \frac{\max_M |f|}{2^*}S_h t_1^{2^*}  < \frac{1}{2} \Phi(t_0) .
 \end{equation}
Let $\vp$ be as in \eqref{eq:preuveexis} with $\| \vp \|_{H^1_h} = 1$. By \eqref{t1} there holds $\Vert t_1 \vp \Vert_{H^1_h} < t_0$. From the definition of $I^{2^*}$ as in \eqref{fonctionnelleq} and by \eqref{eq:preuveexis} there holds:
\begin{equation} \label{ineq3}
\begin{aligned}
I^{2^*}(t_1 \varphi) &  \le  \frac{1}{2}t_1^2 + \frac{\max_M |f|}{2^*}S_h t_1^{2^*} 
 + \frac{1}{2^*} t_1^{-2^*} 
    \frac{C(n)}{(S_h \max_M |f|)^{n-1}} \\
\end{aligned}
\end{equation}
which becomes, with \eqref{C(n)}, \eqref{t1} and \eqref{valt1}:
\begin{equation} \label{mp}
I^{2^*}(t_1 \varphi)  < \frac{1}{2} \Phi(t_0) + \frac{1}{2n} \left( S_h \max_M |f| \right)^{- \frac{n-2}{2}}  < \Phi_0(t_0) .
\end{equation}
By \eqref{limphiq} and with \eqref{ineq2} and \eqref{mp} we can apply Theorem \ref{Th01} and obtain two solutions of \eqref{EL} for $3 \le n \le 5$. 
As one can check, condition \eqref{eq:preuveexis} is satisfied when $\vp$ is a positive constant equal to $ \int_M h dv_g $ whenever $a$ satisfies 
\[\Vert a \Vert_\infty < C(n) \left( S_h \max_M |f| \right)^{1-n} {V_g}^{-1} \left( \int_M h dv_g \right)^{- \frac{2^*}{2}}\hskip.1cm , \] 
where $V_g$ is the volume of $(M,g)$. We have thus just proved that for $\theta$ small enough, equation \eqref{eq:einlicht} has two smooth positive solutions. Thus, the parameter $\theta_1$ appearing in the statement of Theorem \ref{Th1} can be defined as
\begin{equation} \label{theta1}
\theta_1 = \sup ~\left \{~ \xi~ \textrm{ such that for any } \theta \in [0,\xi] \textrm{ there exist at least two solutions of \eqref{eq:einlicht}}\right \}
\end{equation}
and \eqref{eq:preuveexis} provides us with a lower bound on $\theta_1$:
\[ \theta_1 \ge \inf \frac{C(n)}{\left( \max_M |f|  S_h \right)^{n-1}} \left( \int_M \frac{a}{\varphi^{2^*}} dv(g) \right)^{-1}\hskip.1cm , \]
where the infimum is taken over all the smooth positive function $\vp$ on $M$ with $\Vert \vp \Vh = 1$ and $C(n)$ is as in \eqref{C(n)}.

\subsection{The Intermediate case.}  \label{chcroissant}
Just as we defined $\theta_1$ in \eqref{theta1}, we let:
\begin{equation} \label{theta2}
\theta_2 = \sup ~\{~  \theta > 0 \textrm{ such that equation \eqref{eq:einlicht} admits a smooth positive solution } \}.
\end{equation}
There holds $\theta_1 \le \theta_2$ and for any $\theta < \theta_2$ there exists at least one solution of \eqref{eq:einlicht}. Indeed, by definition of $\theta_2$ in \eqref{theta2}, for any $\theta  < \theta_2$ there exists $\theta < \theta_0 \le \theta_2$ such that a solution $u_{\theta_0}$ of \eqref{eq:einlicht} for $\theta = \theta_0$ exists. For any $\xi < \theta_0$, $u_{\theta_0}$ is then a super solution of
\begin{equation} \label{elxi} \tag{$EL_{\xi}$}
\triangle_g u + h u = f u^{2^*-1} + \frac{\xi a}{u^{2^*+1}}.
\end{equation}
Since there exist sub solutions of \eqref{elxi} as small as we want by \eqref{soussol}, \eqref{elxi} has at least one smooth positive solution.

\subsection{The Positive case} \label{caspositif}
It is known that $\theta_2$ is finite when $f > 0$, see Theorem $2.1$ in Hebey-Pacard-Pollack \cite{HePaPo}. To conclude the proof of Theorem \ref{Th1} we now show, when $f >0$, that for every $\theta < \theta_2$ equation \eqref{eq:einlicht} has at least two distinct solutions and that for $\theta = \theta_2$ there is exactly one solution. We define for positive $u \in H^1(M)$ the energy functional of \eqref{eq:einlicht}: 
\begin{equation} \label{Itheta}
I_\theta(u) = \frac{1}{2} \int_M \left( |\nabla u|_g^2 + hu^2 \right) dv_g - \frac{1}{2^*} \int_M f u^{2^*} dv_g + \frac{\theta}{2^*} \int_M \frac{a}{u^{2^*}}dv_g .
\end{equation}
For any $\theta < \theta_2$, Proposition \ref{solmin} provides us with a canonical minimal solution which we shall call $\vp_\theta$. It is stable, that is to say that $D^2 I_\theta(\vp_\theta)$ is a nonnegative bilinear form. Now we investigate the path $\theta \mapsto \vp_\theta$ and show that $D^2 I_\theta(\vp_\theta)$ is actually a definite positive bilinear form if $ \theta < \theta_2$. This will give us a suitable mountain-pass structure to apply Theorem \ref{Th01}. Proposition \ref{solmin} shows that the path $\theta \mapsto \vp_\theta$ is increasing: $\vp_\theta(x) < \vp_\eta(x)$ for all $x \in M$ and $\theta < \eta$. We first prove that $\theta \mapsto \vp_\theta$ is continuous: 

\begin{lemme} \label{c0infini}
The path $\theta \mapsto  \vp_\theta$ is continuous for the $L^\infty(M)$-norm, that is to say:
\begin{equation} \label{eqlim}
 \Vert  \vp_{\theta_0} -  \vp_\theta \Vert_\infty \to 0 
 \end{equation}
as $\theta \to \theta_0$, for any $0 < \theta_0 < \theta_2$, where $\theta_2$ is as in \eqref{theta2}. 
\end{lemme}

\begin{proof}
Let $(\theta_k)_k$ be a sequence of positive numbers converging to some $0 < \theta_0 < \theta_2$; up to assuming $k$ large enough, we can assume that $0 < \theta_0 - \ve \le \theta_k \le \theta_0 + \ve < \theta_2$ for all $k$, for some positive $\ve < (\theta_2 - \theta_0)/2$ . By Proposition \ref{solmin} we have
\begin{equation} \label{bornek}
 \vp_{\theta_0 - \ve} \le \vp_{\theta_k} \le \vp_{\theta_0 + \ve}\hskip.1cm ,
 \end{equation}
where $\vp_{\theta_k}$ denotes the minimal solution of \eqref{eq:einlicht} for $\theta = \theta_k$. By standard elliptic theory and \eqref{bornek} there exists $\vp_0$ in $C^2(M)$ solution of 
\[ \triangle_g \vp_0 + h \vp_0  = f \vp_0^{2^*-1} + \theta_0 a \vp_0^{-2^*-1} \]
such that $\vp_{\theta_k} \to \vp_0$ in $C^2(M)$ as $k \to \infty$. By definition of the minimal solution as given by Proposition \ref{solmin}, there holds $\vp_0 \ge \vp_{\theta_0}$. We proceed by contradiction and assume that $\vp_{\theta_0} < \vp_0$. We then define for any $t \in [0;1]$
\[ m(t) =  I_{\theta_0} \big( t \vp_{\theta_0}  + (1-t) \vp_0 \big)\hskip.1cm , \]
where $I_{\theta_0}$ is as in \eqref{Itheta}. 
By proposition \ref{solmin}, for any $k$, $ \vp_{\theta_k}$ is stable, that is $D^2 I_{\theta_k}(\vp_{\theta_k}) \ge  0$, and thus $\vp_0 $ is stable. Hence $m''(0) \ge 0$. Using \eqref{Itheta} we can compute $m^{(3)}(t)$ for any $t \in [0,1]$, where $m^{(3)}$ is the third derivative of $m$. There holds
\[ \begin{aligned}
m^{(3)}(t)   = -(2^*-1)(2^*-2) \int_M f \big(t \vp_{\theta_0} + (1-t) \vp_0 \big)^{2^*-3}(\vp_{\theta_0} - \vp_0)^3 dv_g \\
- (2^*+1)(2^*+2) \int_M \theta_0 a \big( t \vp_{\theta_0} + (1-t) \vp_0 \big)^{-2^*-3} ( \vp_{\theta_0} - \vp_0)^3 dv_g .\\
\end{aligned} \]
Since $\vp_0 > \vp_{\theta_0}$ and $f > 0$, $m^{(3)}(t)$  is positive for all $t \in (0,1)$. Hence $m''$ is a positive function of $t$ for $0 < t \le 1$ and $m'$ is increasing in $(0,1)$. But this is impossible since both $\vp_0$ and $\vp_{\theta_0}$ are solutions of \eqref{eqlim} and there thus holds $m'(0) = m'(1) = 0$. Hence $\vp_0 = \vp_{\theta_0}$ and $\theta \mapsto \vp_\theta$ is continuous.
\end{proof}

\noindent For any $\theta, \eta < \theta_2$ we have:
\begin{equation} \label{diff}
\begin{aligned}
\triangle_g (\vp_\theta - \vp_{\eta}) + h (\vp_\theta - \vp_{\eta}) = f \left( \vp_\theta^{2^*-1} - \vp_{\eta}^{2^*-1} \right) + a \theta \left( \frac{1}{\vp_\theta^{2^*+1}} - \frac{1}{\vp_{\eta}^{2^*+1}} \right) \\
+ (\theta - \eta) \frac{a}{\vp_\eta^{2^*+1}} .\\
\end{aligned}
\end{equation}
After multiplication of \eqref{diff} by $\vp_\theta - \vp_\eta$ and integration, using Lemma \ref{c0infini} and since $\triangle_g +h$ is coercive we get that $\theta \mapsto \vp_\theta$ is continuous for the $H^1$ norm. We let $\lambda(\theta)$ be the first eigenvalue of the linearized operator at $\vp_\theta$ and $u_\theta$ be the positive associated eigenvector with $\Vert u_\theta \Vert_{\infty} = 1$, thus satisfying:
\begin{equation} \label{valeurpropre}
\triangle_g u_\theta + h u_\theta - (2^*-1)f\vp_\theta^{2^*-2} u_\theta + (2^*+1) \frac{\theta a}{\vp_\theta^{2^*+2}} u_\theta = \lambda(\theta) u_\theta .
\end{equation}
To prove that $D^2 I_\theta(\vp_\theta)$ is positive-definite for any $\theta$ we proceed by contradiction and assume that for some $0 < \theta_0 <\theta_2$ there holds $\lambda(\theta_0) = 0$. We let $\vp_0 = \vp_{\theta_0}$ and $u_0 = u_{\theta_0}$. Our goal is to obtain an asymptotic expansion of $\vp_\theta$ in $H^1(M)$. First, we claim that there holds: 
\begin{equation} \label{quotientH1}
 \frac{\Vert \vp_\theta - \vp_0 \Vert_{H^1_h}}{\theta-\theta_0} \to +\infty 
 \end{equation}
as $\theta \to \theta_0$. Indeed, we define for all $\theta \not = \theta_0$:
\begin{equation} \label{psit}
\psi_\theta = \frac{ \vp_\theta - \vp_0 }{\theta-\theta_0}.
\end{equation}
First, since $\lambda(\theta_0) = 0$, $\theta \mapsto \vp_\theta \in H^1(M)$ is not differentiable at $\theta_0$, otherwise differentiating \eqref{eq:einlicht} with respect to $\theta$ at $\theta_0$ would yield a contradiction with \eqref{valeurpropre}. Hence $\psi_\theta$ has no limit in $H^1(M)$ for $\theta$ going to $\theta_0$. Moreover, $(\psi_\theta)_\theta$ is not even bounded for the $H^1$ norm up to a subsequence. If we assume the contrary, there exists $\psi_0 \in H^1(M)$ such that $\psi_\theta$ converges weakly in $H^1(M)$, up to a subsequence, to $\psi_0$. Then \eqref{diff} gives
\begin{equation} \label{derive}
\begin{aligned}
 \triangle_g \psi_\theta + h \psi_\theta = f \frac{ \vp_\theta^{2^*-1} - \vp_{0}^{2^*-1} }{\vp_\theta - \vp_0} \psi_\theta + a \theta \left( \frac{1}{\vp_\theta^{2^*+1}} - \frac{1}{\vp_{0}^{2^*+1}} \right) (\vp_\theta - \vp_0)^{-1}\psi_\theta + \frac{a}{\vp_0^{2^*+1}}\hskip.1cm .
 \end{aligned}
 \end{equation}
Now we integrate \eqref{derive} against the eigenvector $u_0$ and let $\theta$ go to $\theta_0$: since $\psi_\theta$ converges weakly to $\psi_0$ we can use Lebesgue's dominated convergence theorem to get
\[ \int_M \langle  \nabla \psi_0 , \nabla u_0 \rangle dv_g  + \int_M \left[ h - (2^*-1) f \vp_0^{2^*-2} + (2^*+1) \frac{a \theta_0}{\vp_0^{2^*+2}} \right] \psi_0 u_0 = \int_M \frac{a}{\vp_0^{2^*+1}} u_0 dv_g > 0 \hskip.1cm ,\]
and this yields a contradiction with \eqref{valeurpropre}.  Hence $\Vert \psi_\theta \Vert_{H^1} \to + \infty$ as $\theta \to \theta_0$ and this proves \eqref{quotientH1}. We now come back to \eqref{diff}, pick $\eta = \theta_0$, integrate it against $\vp_\theta - \vp_0$ and use once again Lebesgue's dominated convergence theorem. We obtain for any $\theta < \theta_2$:
\begin{equation} \label{normes}
\begin{aligned}
&\Vert \vp_\theta - \vp_0 \Vert_{H^1_h}^2 \le (\theta - \theta_0)  \int_M \frac{a}{\vp_0^{2^*+1}} (\vp_\theta - \vp_0) dv_g \\
&  + \left( (2^*-1)\int_M f \vp_0^{2^*-2} dv_g + \theta_0 \int_M (2^*+1)\frac{a}{\vp_0^{2^*+2}} dv_g + o_{\theta \to \theta_0}(1)\right) \Vert \vp_\theta - \vp_0 \Vert_\infty^2. \\
\end{aligned}
\end{equation}
With \eqref{quotientH1} and since $f >0$ we thus get:
\begin{equation} \label{quotientL}
\Vert \psi_\theta \Vert_\infty 
\to +\infty\hskip.1cm ,
\end{equation}
where $\psi_\theta$ is as in \eqref{psit}. Define for $\theta \not = \theta_0$
\begin{equation} \label{Phit}
\Phi_\theta = \frac{\vp_\theta - \vp_0}{\Vert \vp_\theta - \vp_0\Vert_{\infty} }.
\end{equation}
Using \eqref{quotientL} in \eqref{normes} we see that $\Phi_\theta$ is $H^1$-bounded and then converges weakly, up to a subsequence, to some function $\Phi_0$ in $H^1(M)$. Since $\Vert \Phi_\theta \Vert_\infty = 1$, using \eqref{diff} and \eqref{quotientL} we obtain that 
\begin{equation} \label{conv1}
\Phi_\theta \rightharpoonup \Phi_0 =  \left\{ 
\begin{aligned}
u_0 & \textrm{ if } \theta > \theta_0 \\
- u_0 & \textrm{ if } \theta < \theta_0 \\
\end{aligned}
\right.
\end{equation}
in $H^1(M)$ as $\theta$ goes to $\theta_0$, where $u_0$ is as in \eqref{valeurpropre}. By \eqref{diff} and \eqref{valeurpropre} we can write:
\begin{equation} \label{derive2}
\begin{aligned}
&\left(\triangle_g + h \right) \left( \Phi_\theta - \Phi_0 \right) = f \left( \frac{\vp_\theta^{2^*-1} - \vp_{\theta_0}^{2^*-1}}{\vp_\theta - \vp_0} \right)  \Phi_\theta - (2^*-1)f\vp_0^{2^*-2} \Phi_0  \\
& + a (\vp_\theta - \vp_0)^{-1} \left( \frac{1}{\vp_\theta^{2^*+1}} - \frac{1}{\vp_{\theta_0}^{2^*+1}} \right) \Phi_\theta + (2^*+1) \frac{\theta_0 a}{\vp_0^{2^*+2}} \Phi_0 + \frac{\theta - \theta_0}{\Vert \vp_\theta - \vp_0 \Vert_\infty} \frac{a}{\vp_\theta^{2^*+1}}.  \\
\end{aligned}
\end{equation}
Integrating \eqref{derive2} against $\Phi_\theta - \Phi_0$, using \eqref{quotientL} and using Lebesgue's dominated convergence theorem shows that $\Phi_\theta$ converges strongly to $\Phi_0$ in $H^1(M)$. This gives the following expansion in $H^1(M)$ for $\vp_\theta$ as $\theta \to \theta_0$:
\begin{equation} \label{dl}
\vp_\theta = \vp_0 + \ve_\theta u_0 + o(\ve_\theta)\hskip.1cm ,
\end{equation}
where we have set $\ve_\theta = \mathrm{sgn}(\theta - \theta_0) \Vert \vp_\theta - \vp_0 \Vert_{\infty}$. To conclude it remains to integrate \eqref{valeurpropre} against $u_0$ and use \eqref{dl} and \eqref{quotientL} to obtain, at the first order in $\ve_\theta$:
\begin{equation} \label{bigth9}
\begin{aligned}
- \left( \int_M \left[ (2^*-1)(2^*-2) f \vp_0^{2^*-3} u_0^3 + (2^*+1)(2^*+2) \frac{\theta_0 a}{\vp_0^{2^*+3}} u_0^3 \right] dv_g \right) \ve_\theta + o(\ve_\theta) \\
  = (1 + o(1)) \lambda(\theta) \int_M u_0^2 dv_g. \\
\end{aligned}
\end{equation}
As one can check, to obtain \eqref{bigth9} we used that $u_\theta \to u_0$ in $L^2(M)$ as $\theta \to \theta_0$. To prove this, first recall the variational characterization of $\lambda(\theta)$ for all $\theta$ as:
\[ \lambda(\theta) = \inf_{\psi \in H^1(M), \| \psi \|_2 = 1} \int_M \left( |\nabla \psi|^2 + \left[ h - (2^*-1) f \vp_\theta^{2^*-2} + (2^*+1) \frac{\theta a}{\vp_\theta^{2^*+2}}\right] \psi^2\right) dv_g. \]
In particular, this shows that $\lambda: \theta \mapsto \lambda(\theta)$ is upper semi-continuous. Since by its very definition $\vp_\theta$ is stable, $\lambda(\theta)$ is always non negative and there thus holds $\lambda(\theta) \to 0$ as $\theta \to \theta_0$. Then the convergence of $u_\theta$ towards $u_0$ in $L^2(M)$ holds by standard elliptic theory since $u_\theta$ satisfies in addition \eqref{valeurpropre} and $\| u_\theta \|_\infty = 1$ for all $0 < \theta < \theta_2$. 
Since we assumed $f>0$ we see with \eqref{bigth9} that $\lambda$ changes sign at $\theta_0$. By assumption $\lambda(\theta)$ is always non negative: we have a contradiction and $D^2 I_\theta(\vp_\theta)$ is definite positive for all $0 < \theta < \theta_2$.
Application of Theorem \ref{Th01} is then straightforward and gives a second solution of \eqref{eq:einlicht} when $0 < \theta < \theta_2$. Letting $\theta_\star = \theta_2$ , we have thus shown that, if $f >0$, \eqref{EL} has at least two solutions for $\theta < \theta_\star$ and none for $\theta > \theta_\star$. 
Now we show that in the limit case $\theta = \theta_\star$ there is a unique solution. Its existence is obtained with a stability argument and its uniqueness is a consequence of the non coercivity of $D^2 I_{\theta_\star}(\vp_{\theta_\star})$, where $I_{\theta_\star}$ is as in \eqref{Itheta}. 
First, a solution exists for $\theta = \theta_\star$. Indeed, let $(\theta_k)_k$, $0 < \theta_k < \theta_\star$ be a sequence converging to $\theta_\star$ (which is finite since $f > 0$). Then by Theorem \ref{Thstabi} applied with $q_k = 2^*$ and $a_k = \theta_k a$ in the case $f > 0$, the sequence $\vp_{\theta_k}$ converges in $C^{1, \alpha}(M)$ to some smooth positive function $\vp_{\theta_\star}$ that solves  \eqref{eq:einlicht} with $\theta = \theta_\star$. To prove uniqueness, we first show that the minimal solution $\vp_{\theta_\star}$ is not strictly stable. For the sake of simplicity, we let  $\vp_\star = \vp_{\theta_\star}$. We proceed by contradiction and assume that $\vp_\star$ is strictly stable: that is,
we assume that 
\begin{equation} \label{coercif}
\lambda(\theta_\star) > 0\hskip.1cm ,
\end{equation}
where $\lambda(\theta_\star)$ is the first eigenvalue of the linearization of \eqref{ELt} for $\theta = \theta_\star$ at $\vp_\star$. 
Let $0 < \alpha < 1$: since we assumed that $\triangle_g + h$ is coercive it is, by standard elliptic theory, an isomorphism between $C^{k+2,\alpha}(M)$ and $C^{k, \alpha}(M)$ for any integer $k \ge 0$. Let $ \delta = \frac{1}{2} \inf_M \vp_\star >0$, $B(\vp_\star, \delta) = \{ v \in C^{2,\alpha}(M) \hbox{ s.t. } \| v-\vp_\star\|_{C^{2,\alpha}} \le \delta \}$ and $\Omega =  B(\vp_\star, \delta) \times \mathbb{R} $. The following mapping is thus well-defined:
\begin{equation} \label{appliF}
F: \left \{
\begin{aligned}
\Omega 
& \longrightarrow  C^{2,\alpha}(M) \\ (u,\theta) & \longmapsto  u - \left( \triangle_g + h \right)^{-1} \left(f u^{2^*-1} + \frac{\theta a}{u^{2^*+1}} \right).
\end{aligned} \right.
\end{equation}
It is easily seen that $F$ is $C^1$ in $ \Omega$ and that there holds, for all $(u,\theta) \in \Omega$ and for all $v \in C^{2, \alpha}(M)$:
\[ \frac{\partial F}{\partial u}(u,\theta)(v) = v - \left( \triangle_g + h \right)^{-1} \left( (2^*-1) f u^{2^*-2} v - (2^*+1) \frac{\theta a}{u^{2^*+2}} v \right) .\]
We claim that $ \frac{\partial F}{\partial u}(\vp_\star,\theta_\star)$ is an isomorphism. Indeed, let $w \in C^{2,\alpha}(M)$. Finding $v\in C^{2,\alpha}(M)$ such that 
\begin{equation} \label{bij}
 \frac{\partial F}{\partial u}(\vp_\star,\theta_\star)(v) = w 
 \end{equation}
is equivalent to solve
\begin{equation} \label{isomorphisme}
 \triangle_g v+ \left[h - (2^*-1) f \vp_\star^{2^*-2} + (2^*+1) \frac{\theta_\star a}{\vp_\star^{2^*+2}} \right] v = (\triangle_g + h) w .
 \end{equation}
By \eqref{coercif} the left-hand side is a coercive linear operator. Since $\triangle_g w+ hw \in C^{0,\alpha}(M)$ there exists, by standard elliptic theory, a unique $v \in C^{2,\alpha}(M)$ satisfying \eqref{isomorphisme}, and hence \eqref{bij}. As one can check, $F(\vp_\star, \theta_\star) = 0 $ and more generally, for any $(u,\theta) \in \Omega$, $u$ solves \eqref{ELt} if and only if $F(u,\theta) = 0$. With \eqref{bij} the implicit function theorem applies to $F$ at $(\vp_\star,\theta_\star)$ and shows that the set of solutions of the equation $F(u,\theta) = 0$ is locally a path in $B(\vp_\star, \delta)$ parameterized by $\theta \in (\theta_\star- \ve, \theta_\star + \ve)$ for some positive $\ve$. But this is impossible since no solutions of \eqref{ELt} exist if $\theta > \theta_\star$. 

\medskip

\noindent Hence $\lambda(\theta_\star) = 0$ and there exists a positive $\eta \in C^\infty(M)$ such that
\begin{equation} \label{1vp}
 \triangle_g \eta+ \left[h - (2^*-1) f \vp_\star^{2^*-2} + (2^*+1) \frac{\theta_\star a}{\vp_\star^{2^*+2}} \right] \eta = 0.
\end{equation}
The uniqueness of solutions for $\theta = \theta_\star$ then follows as in Ma-Wei \cite{MaWei}. If there exists a solution $v$ different from $\vp_\star$, by minimality of $\vp_\star$ there would hold $v > \vp_\star$. Since $f >0$ in $M$ and $a \ge 0$ in $M$, for any $x\in M$ the functions $u \in \RR \mapsto f(x) u^{2^*-1}$ and $u \in \RR \mapsto \frac{\theta_\star a(x)}{u^{2^*+1}}$ are respectively strictly convex and convex. We therefore obtain:
\begin{equation} \label{signe}
\begin{aligned}
\triangle_g (v - \vp_\star) + \left[ h - (2^*-1) f \vp_\star^{2^*-2} + (2^*+1) \frac{\theta_\star a}{\vp_\star^{2^*+2}} \right](v - \vp_\star) > 0 
\end{aligned}
\end{equation}
and integrating \eqref{signe} against $\eta$ and using \eqref{1vp} gives a contradiction. This concludes the proof of Theorem \ref{Th1} . 

\bibliographystyle{abbrv}
\bibliography{biblio}

\end{document}